\documentclass[11pt]{article}
\usepackage{ulem}
\usepackage[makeroom]{cancel}
\usepackage{color}
\usepackage{xcolor}
\usepackage{amsfonts}
\usepackage{amssymb}
\usepackage{amscd}
\usepackage{amsmath}
\usepackage{epsfig}
\usepackage{graphicx, subfigure}
\usepackage{latexsym}
\usepackage{pst-3dplot}
\usepackage{verbatim}
\usepackage{tikz}
\usepackage{pdfsync}
\usepackage{hyperref}
\usepackage{multirow}
\usepackage[utf8]{inputenc}

\usepackage{color}
\usepackage{graphicx}

\usepackage{caption}
\usepackage{lipsum}
\newlength{\dinwidth}
\newlength{\dinmargin}
\newtheorem{definition}{Definition}
\newtheorem{theorem}{Theorem}

\newtheorem{proposition}{Proposition}
\newtheorem{corollary}{Corollary}
\newtheorem{remark}{Remark}
\newtheorem{lemma}{Lemma}
\newtheorem{example}{Example}

\def \surf{{\cal L}}

\def \i{{\rm i}}
\def\Axi{A_{x_i}}
\def\Aum{A_{u_m}}

\def\Aaj{A_{a_j}}

\def\pxi{P_{x_i}}
\def\pxn{P_{x_n}}
\def\pxk{P_{x_k}}
\def\pxj{P_{x_j}}
\def\pum{P_{u_m}}
\def\pualpha{P_{u_\alpha}}
\def\puk{P_{u_k}}
\def\puj{P_{u_j}}
\def\paj{P_{a_j}}
\def\pak{P_{a_k}}

\def\N{M_2}
\def\M{M_1}

\def\pualpha{P_{u_\alpha}}
\def\x{{\bf x}}
\def\u{{\bf u}}
\def\curve{ T}
\def\surf{\mathcal T}
\def\proj{\hat{ T}}

\def\Tx{\surf_\x}
\def\Txu{\surf_{\x\u}}
\def\T{\mathcal T}

\def\l{\lambda}

\def\a{{\bf a}}
\def\b{{\bf b}}
\def\x{{\bf x}}

\def\B{\mathcal B}

\makeatletter
\@addtoreset{equation}{section}
\makeatother

\begin{document}
\title{Isoperiodic deformations of Toda curves and chains, the difference Korteweg - de Vries equation, and $SU(N)$ Seiberg-Witten theories}

\author{Vladimir Dragovi\'c$^1$ and Vasilisa Shramchenko$^2$}
\date{}

\maketitle

\footnotetext[1]{Department of Mathematical Sciences, The University
	of Texas at Dallas, 800 West Campbell Road, Richardson TX 75080,
	USA. Mathematical Institute SANU, Kneza Mihaila 36, 11000
	Belgrade, Serbia.  E-mail: {\tt
		Vladimir.Dragovic@utdallas.edu}--the corresponding author}

\footnotetext[2]{Department of mathematics, University of
	Sherbrooke, 2500, boul. de l'Universit\'e,  J1K 2R1 Sherbrooke, Quebec, Canada. E-mail: {\tt Vasilisa.Shramchenko@Usherbrooke.ca}}

\

\

{\it  Dedicated to the memory of Igor Moiseevich Krichever (1950-2022).}

\

\begin{abstract} We introduce the dynamics of Toda curves of order $N$ and derive  differential  equations governing this dynamics.
We apply the obtained results to describe isoperiodic deformations of $N$-periodic Toda  chains and periodic difference Korteweg-de Vries equation.  We describe deformations of the essential spectra of $N$-periodic two-sided Jacobi matrices. We also study singular regimes  of $SU(N)$ Seiberg-Witten theory and describe their deformations preserving the  number  of singularities where new massless particles may occur. We introduce and describe isoequilibrium deformations of arbitrary collections of $d$ real disjoint closed intervals. We conclude by providing explicit triangular solutions to constrained Schlesinger systems.
	
	\vskip 1cm
	
	
	
\end{abstract}

\

\tableofcontents

\

\section{Introduction}
\label{sect_introduction}

    Toda curves of order $N$ are real hyperelliptic curves, as well as the corresponding compact Riemann surfaces, having two points $\infty^-$ and $\infty^+$ at infinity, carrying a meromorphic function with a pole of order $N$ at $\infty^-$ and a zero of order $N$ at $\infty^+$ and no other zeros or poles \cite{HillToda}.  Toda curves encode the equations of the periodic Toda chain and give rise to solutions of those equations \cite{Krich1981}. Toda curves  arise also from polynomial solutions to the Pell equation, some of which are the generalized Chebyshev polynomials. The question we address here is finding continuous $g$-parameter families of Toda curves in the moduli space of all hyperelliptic curves. This question can be reformulated in terms of an Abelian differential of the third kind $\Omega$ having simple poles at $\infty^-$ and $\infty^+$ with residues $-1$ and $+1$, respectively, and no other poles. The families of Toda curves  we describe here are the families of hyperelliptic real surfaces along which the periods of $\Omega$ with respect to some chosen homology basis stay constant. We call these families {\it isoperiodic families relative to $\Omega$}.
   \

Denoting the set of branch points of the hyperelliptic coverings of $\mathbb CP^1$ by $\{0,1,x_1, \dots, x_g, u_1, \dots, u_g\}$, we show that isoperiodic families can be parametrized by $g$ branch points. We derive a system of differential equations of order two for $u_1, \dots, u_g$ as functions of $x_1, \dots, x_g$ that are satisfied along isoperiodic families. Notably, these equations are rational in $u_j, x_j$ and polynomial in the first order derivatives of $u_j$.  We also prove that, under some genericity assumptions, the families of hyperelliptic curves with above branch points for which $u_j(x_1, \dots, x_g)$ for $j=1, \dots, g$ satisfy our system of differential equations  are isoperiodic.
In the case of families of genus one curves we obtain an ordinary differential equation of the second order with rational coefficients, see \eqref{iso_genus1}.

The existence of isoperiodic families of hyperelliptic surfaces implies the existence of continuous families of Toda curves and of solutions to the Pell equation. On the other hand, periods of meromorphic differentials of the third kind appear in algebro-geometric solutions of equations of mathematical physics, such as  the Toda chain equations and the difference Korteweg - de Vries equation. Some conditions for the periods of $\Omega$ single out periodic solutions of these equations. It is difficult though to find Riemann surfaces for which these conditions are satisfied. We show the existence of continuous families of Riemann surfaces with given fixed periods of $\Omega$. Once one such hyperelliptic Riemann surface is known, it suffices to evolve it according to our equations to obtain a continuous family of hyperelliptic surfaces satisfying the required condition for the periods of $\Omega$.

Remarkably, there is a relationship between  hyperelliptic curves of the form
$$
\mu^2=\mathcal P_N^2(z)-1,
$$
where $\mathcal P_N$ is a polynomial of degree $N$, and the Seiberg-Witten theory, see \cite{AF1995, KLTY,GKMMM}.  We extend it here to a fruitful relationship between the Seiberg-Witten theory and Pell's equations. In this case, isoperiodic families of hyperelliptic Riemann surfaces, relative to a differential of the third kind, lead to local continuous deformations  of singular regimes in the $SU(N)$ super-symmetric Yang-Mills theory. In the singular regimes, new  massless particles appear. Our deformations preserve the number of  these singularities where new massless particles in the theory may occur.

\section{Toda curves and Pell equations}
\label{sect_TodaPell}

The Pell equation is the following equation for three polynomials: ${\mathcal P}_N$ of degree $N$, ${\mathcal Q}_{N-g-1}$ of degree $N-g-1$, and  $\Delta_{2g+2}$ a monic polynomial of degree $2g+2$  vanishing at the ends of $g+1$ real intervals
%
%
%
\begin{equation}
\label{Pell}
{\mathcal P}_N^2(u)-\Delta_{2g+2}(u){\mathcal Q}_{N-g-1}^2(u)=1.
\end{equation}
Here $N$ and $g$ are positive integers such that $g\geq 1$ and $N\geq g+1.$
We also say that  $\mathcal P_N$ is a solution of the Pell equation, if there exist polynomials $\Delta_{2g+2}$ and  $\mathcal Q_{N-g-1}$ such that \eqref{Pell} holds.
 In the case when all zeros of $\Delta_{2g+2}$ and of ${\mathcal Q}_{N-g-1}$ are real, solutions $\mathcal P_N$ of the Pell equation are the (generalized) Chebyshev polynomials \cite{SY1992}. They are also called  {\it extremal}  polynomials on the set of $g+1$ real intervals, each connecting two consecutive zeros of $\Delta_{2g+2}$ and for which $|\mathcal P_N|\leqslant 1$. Such a set of $g+1$ real intervals is called {\it support of the Pell equation}.  The extremal polynomials are minimal in the sense of the  following extremality  condition    \cite{Bogatyrev2012} that they satisfy: {\it The monic polynomial obtained from $\mathcal P_N$ by dividing out the leading coefficient is the polynomial that deviates from zero over the given set of intervals less than any other monic polynomial of degree $N$ having the same {\it signature} as $\mathcal P_N$}. The {\it  signature} of $\mathcal P_N$ is a vector with $g+1$ integer components where the $j$-th component gives the number of local maxima and minima of $\mathcal P_N$ inside the $j$-th interval. The support of the corresponding Pell equation is then the maximal subset of the real line for which $\mathcal P_N$ is minimal in the above sense. We use the term {\it support of the  Chebyshev polynomial $\mathcal P_N$} to mean the support of the corresponding Pell equation.

Supports obtained from each other by an affine transformation preserving the orientation of the real line are considered equivalent. Therefore we may assume the leftmost interval of the support to be $[0,1].$
Thus, we consider the polynomial $\Delta_{2g+2}$ from \eqref{Pell} to be
\begin{equation}
\label{delta}
\Delta_{2g+2}(u)=u(u-1)\prod_{j=1}^{g}(u- x_j)(u- u_j),
\end{equation}
and the support of the Pell equation to be of the form $[0,1]\cup(\cup_{g=1}^g[x_j, u_j]).$

 The study of extremal  polynomials was initiated by Chebyshev for the case of one interval. Zolotarev studied minimal polynomials on two intervals. Important results of the theory of generalized Chebyshev polynomials were obtained by Markov, Bernstein, Borel, Akhiezer in the early 20th century. Research in this field continues nowadays with various applications, see eg. \cite{AhiezerAPPROX},  \cite{Akh4},  \cite{Bogatyrev2012}, \cite{PS1999},  \cite{Si2011}, \cite{Si2015}, \cite{SY1992}, \cite{Widom} and references therein.

Associated with each solution of the Pell equation is a hyperelliptic algebraic curve
\begin{equation}
\label{T}
\curve=\{(u,v)\in\mathbb C^2 \; | \; v^2=\Delta_{2g+2}(u)\},
\end{equation}
where $\Delta_{2g+2}(u)$ is given by \eqref{delta} with $x_j$ and $u_j$ fixed real numbers. 
The curve $\curve$ is a Toda curve according to  Definition \ref{def_Toda} given below, see also \cite{HillToda}.

The compact projective curve $\proj\subset \mathbb CP^2$ obtained from $\curve$ by adding a point at infinity $[u:v:\lambda]=[0:1:0]$ is singular at the point at infinity. This singularity is a double point giving rise to two points $\infty^+$ and $\infty^-$ on the  associated compact hyperelliptic Riemann surface $\surf$. This surface  is of genus $g$ and is such that there is a holomorphic map $\pi:\surf\to \mathbb CP^2$ where $\pi(\surf)=\proj$ and
%
\begin{equation}
\label{pi}
\pi:\surf\setminus\pi^{-1}([0:1:0])\to\proj\setminus\{[0:1:0]\}=T
\end{equation}
is a biholomorphism with $\pi^{-1}([0:1:0])=\{\infty^+, \infty^-\}$. Using this map, and given that the algebraic curve \eqref{T} is non-singular, we identify the points of the Riemann surface $\surf$ with those of the projective curve $\proj$ or of algebraic curve $\curve$ as follows. We say that $P=(u, v)$ is a point of the surface $\surf$ if $P=\pi^{-1}([u:v:1]).$

The function $u:\curve\to\mathbb C$ defined naturally on the algebraic curve as the projection to the first coordinate may be extended to the Riemann surface $u:\surf\to \mathbb CP^1$ by setting $u(\infty^{\pm})=\infty.$ This allows us to represent the hyperelliptic Riemann surface  as the surface of the two-fold ramified covering $(\surf, u)$ ramified over the points of the set
\begin{equation}
\label{branch}
\B:=\{0,1,x_1,u_1, \dots, x_g, u_g\},
\end{equation}
called the branch points of the covering. We denote the corresponding ramification points by $P_{a_j}=(a_j, 0)$ for $a_j\in\B$ and $j=1, \dots, 2g+2.$ The two points at infinity $\infty^{\pm}$ are characterized by the condition $v\simeq \pm u^{g+1}$ for $P=(u,v)\simeq \infty^{\pm}.$
\begin{definition}
\label{def_Toda}
A hyperelliptic Riemann surface $\surf$  corresponding to the algebraic curve $\curve $ defined by \eqref{T}, \eqref{delta} with $1<x_j<u_j$ for $1\leq j\leq g$  is called a Toda curve of   order $N$ if there is a nonzero integer $N$ such that $N(\infty^+-\infty^-)$ is a principal divisor.
\end{definition}
Let us now show that a surface defined by \eqref{T}, \eqref{delta} is a Toda curve of  order $N$ if $[0,1]\cup(\cup_{g=1}^g[x_j, u_j])$ is a support of the Pell equation. Let $\mathcal P_N$ be the Chebyshev polynomial of degree $N$ satisfying the Pell equation \eqref{Pell} with $\Delta_{2g+2}$ given by \eqref{delta}. On the surface $\surf$ of the algebraic curve \eqref{T}, consider a function given by  $f(P)=\mathcal P_N(u)+v{\mathcal Q}_{N-g-1}(u)$ for $P=(u,v)\in\surf.$
Denoting by $P^*$ the hyperelliptic involution of the point $P=(u,v)$, that is $P^*=(u,-v)$, we may write the Pell equation as $f(P^*)f(P)=1$  since $f(P^*)=\mathcal P_N(u)-v{\mathcal Q}_{N-g-1}(u).$ Given that the function $f$ has a pole  of order $N$ at $\infty^+$, the property $f(P^*)f(P)=1$ allows us to conclude that $f$ has a zero of order $N$ at $\infty^-=(\infty^+)^*.$ Since $f$ does not have other poles on $\surf$, we conclude that its divisor is $(f)=N(\infty^+-\infty^-)$ and thus $\surf$ is a Toda curve of order $N$.

The converse question, under what condition a Toda curve corresponds to a generalized Chebyshev polynomial, is discussed in Corollary \ref{cor:M20}.
\\

Let us denote by $\{\a_1,\dots, \a_g; \b_1, \dots, \b_g\}$ a fixed canonical homology basis on $\surf,$ with the following intersection indices: $\a_k\circ \a_j=\b_k\circ \b_j=0$ and $\a_k\circ \b_j=\delta_{kj}$.

\subsection{Toda families of Riemann surfaces}
\label{sect_family}

We express the Riemann surface structure on a Toda curve by means of the following {\it standard local charts} on the surface $\surf$  induced by the ramified covering $u:\surf\to \mathbb CP^1:$
\begin{align}
\label{coordinates}
&\zeta_{a_k}(P)=\sqrt{u(P) - a_k} \quad\mbox{if}\quad P\sim P_{a_k}, \quad\mbox{with } a_k\in B,
\nonumber
\\
& \zeta_\infty(P)= \frac{1}{{u(P)}}\quad\mbox{if}\quad P\sim \infty^+ \mbox{ or } P\sim \infty^-,
\\
& \zeta_Q(P)=u(P)-u(Q) \quad\mbox{if}\quad P\sim Q \mbox{ and } Q\notin (\{P_0, \dots, P_{a_{2g+2}}\}\cup\{\infty^+\}+\cup\{\infty^-\}).
\nonumber
\end{align}
We may now consider a variation of the branch points $a_k$, which induces a variation of the curve  \eqref{T} and of the associated Riemann surface. We thus obtain a family of compact Riemann surfaces associated with the family of algebraic curves \eqref{T}. However, these curves may not  be  Toda curves.

In order to construct a continuous family of Toda curves of a given order $N$, we use the correspondence between supports of the Pell equation and Toda curves as well as  Theorems 2.7 and 2.12 from \cite{PS1999}. These theorems imply that there exists a continuous family of supports of the Pell equation of the form $[0,1]\cup(\cup_{g=1}^g[x_j, u_j])$, where the endpoints of the intervals denoted by  $u_j, \;j=1, \dots, g$ are smooth functions of the endpoints denoted by $x_j, \;j=1, \dots, g$, while $(x_1, \dots, x_g)$ varies in some neighbourhood in $\mathbb R^n.$ Moreover, each support $[0,1]\cup(\cup_{g=1}^g[x_j, u_j])$ from such a family corresponds to the Chebyshev polynomial of a fixed degree $N$ and a fixed signature, fixed for all supports of the family. Since any support of the Pell equation corresponds to a Toda curve, we have a continuous family of Toda curves of order $N$, parameterized by the $g$ real parameters $x_1, \dots, x_g.$ Let us introduce notation $\x:=(x_1, \dots, x_g)$ and denote by $\Tx$ the Toda curve (that is the corresponding Riemann surface) from the family corresponding to the values $\x$ of the parameters.

We may assume that the variation of $x_1, \dots, x_g$ is such that a canonical homology basis
 can be chosen on all the curves from the family in such a way that the images of the cycles $\{\a_1,\dots, \a_g; \b_1, \dots, \b_g\}$ by $u$ are fixed for all curves from the family, that is do not depend on $x_1, \dots, x_g$.

 The family of curves $\Tx$ \eqref{T}  admits two sections at infinity, we denote them $s_{\infty^+}$ and $s_{\infty^-}$.   We have $\mu\sim u^{g+1}$ locally near $s_{\infty^+}$ and $\mu\sim -u^{g+1}$ at $s_{\infty^-}$.
 The divisor $N(s_{\infty^+}(\x)-s_{\infty^-}(\x))$ is principal on the corresponding curve $\Tx$ from the family.

More generally, {\it abandoning the requirement of reality of the branch points}, we define a Toda family of which the constructed family of Toda curves is an example.
\begin{definition}
\label{def_Todafamily}
A triple $(T, s_{\infty^+}, s_{\infty^-})$ is called a Toda family of  order $N$, if $ T: \T\to  X$ is a smooth fibration with fibers given by compact hyperelliptic genus $g$ Riemann surfaces  $\Tx$ for $\x\in X$ such that
\begin{itemize}
\item $X\subset\mathbb C^g\setminus \Delta$ with $\Delta=\{(x_1, \dots, x_g)\in\mathbb C^g\;|\; \exists\, k\neq j \mbox{ with } x_k=x_j \}$,
\item $u:\Tx\to\mathbb CP^1$ is a function of degree two on $\Tx,$
\item $\{\a_1,\dots, \a_g; \b_1, \dots, \b_g\}$ is a canonical homology basis on $\Tx$ for $\x\in X$ where $u(\a_k)$ and $u(\b_k)$ are independent of $\x\in X$ for all $1\leq k\leq g,$
\item $s_{\infty^+}$ and $s_{\infty^-}$ are two sections of $T$ such that $N(s_{\infty^+}(\x)-s_{\infty^-}(\x))$ is a principal divisor on  $\Tx.$
\end{itemize}
\end{definition}

In the rest of the paper we assume that the Riemann surfaces from a Toda family correspond to hyperelliptic curves of the equation \eqref{T}, \eqref{delta} with $\x=(x_1, \dots, x_g)$ varying in some appropriate subset $X$ of $\mathbb C^g.$ In Section \ref{sect_isoeqs} we derive a system of differential equations for $u_1, \dots, u_g$ as functions of $x_1, \dots, x_g$. All subsequent results for Toda families are valid in particular for Toda families with real branch points corresponding to supports of the Pell equation.

\section{Meromorphic differentials on Toda curves}

\label{sect_differentials}

Let $\Tx$ be a Riemann surface of genus $g$ from a Toda family corresponding to the hyperelliptic curve of the equation
\begin{equation}
\label{surf}
v^2=u(u-1)\prod_{j=1}^{g}(u- x_j)(u- u_j),
\end{equation}
where $\x=(x_1, \dots, x_g), (u_1, \dots, u_g)\in\mathbb C^g.$
 Let us denote $\omega=(\omega_1, \dots, \omega_g)^T$  the vector of holomorphic differentials normalized by the condition $\oint_{\a_j}\omega_k=\delta_{jk}$ with respect to $a$-cycles of the chosen canonical homology basis of the family as in Definition \ref{def_Todafamily}, $j,k=1, \dots, g$.  The {\it $b$-periods} of these differentials define the {\it Riemann matrix} $\mathbb B$ of $\Tx:$
\begin{equation}
\label{B}
\oint_{\b_j}\omega_k=\mathbb B_{jk}, \qquad j,k=1, \dots, g,
\end{equation}
which is symmetric and has a strictly positive imaginary part.

Note that, due to the Abel theorem,  Definition \ref{def_Toda} of a Toda curve  of order $N$  implies  that there exist  rational column vectors $\M$ and $\N$ such that $N\M, N\N\in\mathbb Z^g$ and
\begin{equation}
\label{Abel}
\int_{\infty^-}^{\infty^+}\omega = \M+\mathbb B \N\,.
\end{equation}
Another basis in the space of holomorphic differentials on $\Tx$ consists of differentials $\phi, u\phi, \dots, u^{g-1}\phi$ where
\begin{equation}
\label{phi}
\phi(P)=\frac{du(P)}{v}.
\end{equation}
For our purposes, it is convenient to define a third basis in the space of holomorphic differentials on $\Tx$ formed by  differentials $v_1, \dots, v_g$ defined by
\begin{equation}
\label{v}
v_k(P)=\frac{\varphi (P)\prod_{\beta=1, \beta\neq k}^{g} (u-u_\beta)}{\varphi (P_{u_k})\prod_{\alpha=1, \alpha\neq k}^{g} (u_k-u_\alpha)}, \qquad   k=1, \cdots , g\,.
\end{equation}
These explicit expressions show that the zeros of $v_k$ at the ramification points $\pum$ with $m\neq k$ are of second order and thus there are no other zeros; in particular, the differentials $v_k$ do not vanish at the points at infinity $\infty^+$ and $\infty^-$.

A convenient tool for working with meromorphic differentials is the so-called Riemann fundamental bidifferential $W(P,Q)$, which can be defined on surfaces $\Tx$ of the Toda family for the given choice of a canonical homology basis $\{\a_1,\dots, \a_g; \b_1, \dots, \b_g\}$ as in Definition \ref{def_Todafamily}, and two points $P$ and $Q$ on $\Tx$. The bidifferential can be expressed in terms of an odd theta-function with non-degenerate characteristics or defined by its three properties as follows:
\begin{itemize}
\item Symmetry: $W(P,Q) = W(Q,P);$
\item A second order pole along the diagonal $P=Q$ with the following local expansion in terms of a local parameter $\zeta$ near $P=Q$ and no other singularities:
\begin{equation}
\label{W}
W(P,Q) \underset{P\sim Q}{=} \left( \frac{1}{(\zeta(P) - \zeta(Q))^2}  + {\cal O}(1) \right)d\zeta(P) d\zeta(Q);
\end{equation}
\item Normalization by vanishing of the $a$-periods: $\oint_{\a_j} W(P,Q) = 0$ for all $j=1, \dots, g.$ Due to the symmetry,  these integrals can be computed with respect to either $P$ or $Q$.
\end{itemize}
This definition can be shown to imply
\begin{equation}
\label{bW}
\oint_{\b_j}W(P,Q) = 2\pi{\rm i}\,\omega_j(P), \qquad j=1, \dots, g,
\end{equation}
for the holomorphic normalized differentials $\omega=(\omega_1, \dots, \omega_g)^T.$

In this paper, the main role is played by the meromorphic differential of the third kind defined by
\begin{equation}
\label{Omega}
\Omega(P) = \Omega_{\infty^+\infty^-}(P) -2\pi\i \N^T\omega,
\end{equation}
where $\N\in\mathbb Z^g$ is a constant column vector from \eqref{Abel} and $\Omega_{\infty^+\infty^-}$ is the  differential of the third kind having simple poles at $\infty^+$ and $\infty^-$ with residues $+1$ and $-1$, respectively, and normalized by vanishing of its $a$-periods, $\oint_{\a_k}\Omega_{\infty^-\infty^+}$ for all $k=1, \dots, g$. This latter differential can be defined through the fundamental bidifferential as follows
\begin{equation}
\label{OW}
\Omega_{\infty^-\infty^+}(P) = \int_{\infty^-}^{\infty^+}W(P,Q) \,.
\end{equation}
Note that the $b$-periods of \eqref{OW} are
%
$\oint_{b_k}\Omega_{\infty^-\infty^+}=2\pi\i\int_{\infty^-}^{\infty^+}\omega_k$
%
and thus due to \eqref{Abel}, the $b$-periods of $\Omega$ \eqref{Omega} are given by
\begin{equation}
\label{bOmega}
\oint_{\b_k}\Omega =2\pi\i {\M}_k\,,
\end{equation}
where ${\M}_k$ is the $k$-th component of $\M$.
The $a$-periods of $\Omega$ are constant and equal to the components of the vector $\N$.

We need to work with meromorphic differentials {\it evaluated} at ramification points of the covering $u:\surf\to \mathbb CP^1.$  The evaluation is defined through the standard local parameters from Section \ref{sect_family} as follows. Let  $\Upsilon$ be a meromorphic differential on $\surf\,.$ We define its {\it value} at a point $ Q_o\in\surf$ which is not a pole of $\Upsilon$ as the constant term of the Taylor series expansion of the differential with respect to the standard local parameter $\zeta$ from the list \eqref{coordinates} at $ Q_o$, that is
\begin{equation}
\label{evaluation}
\Upsilon( Q_o) = \frac{\Upsilon(P)}{d\zeta(P)}\Big{|}_{P= Q_o}.
\end{equation}
As an example, for the value at a ramification point $P_{a_j}$ with $a_j\in B$ of the holomorphic differential \eqref{phi}, we use the local parameter $\zeta_{a_j}(P) = \sqrt{u(P)-a_j}$ for which we have $du(P)=2\zeta_{a_j}(P)d\zeta_{a_j}(P)$.
Thus
\begin{equation}
\label{phi-eval}
\phi(\paj) =\frac{\phi(P)}{d\zeta_{a_j}(P)}\Big{|}_{\zeta_{a_j}=0}= \frac{2}{\sqrt{\prod_{\substack{i=1\\i\neq j}}^{2g+1}(a_j-a_i)}}.
\end{equation}
Using this evaluation, we can characterize  differentials $v_1, \dots, v_g$ \eqref{v} as the holomorphic differentials normalized by the values at $g$ points $P_{u_j}$ on the surface as follows:
\begin{equation}
\label{vcond}
v_k(P_{u_j}) =\delta _{kj},   \quad\text{with} \;\;  1\leq k,j \leq g.
\end{equation}
This definition of evaluation  extends to the bidifferential $W$ by considering one of the arguments of $W$ fixed and thus obtaining a differential with respect to the other argument. Thus,
\begin{equation}
\label{W-eval}
W(P, Q_o) =\frac{W(P, Q)}{d\zeta_{Q_o}(Q)}\Big{|}_{\zeta_{Q_o}=0},
\end{equation}
with $\zeta_{Q_o}$ being the standard local parameter in a neighbourhood of $Q_o$ defined in Section \ref{sect_family}.

We end this section by obtaining some useful identities for the introduced differentials.
\begin{lemma} Let $\Tx$ be the compact Riemann surface corresponding to the algebraic curve \eqref{surf} and let $P_{a_j}$ denote the ramification point of the covering $u:\Tx\to \mathbb CP^1$  corresponding to the branch point $a_j\in \B=\{0,1,x_1, u_1, \dots, x_g, u_g\}.$ Let the differentials $\Omega$, $\phi$, $v_j$ be the differentials defined respectively by \eqref{Omega}, \eqref{phi}, \eqref{v} on $\Tx$ and $W$ be the Riemann bidifferential on $\Tx$.  Then  for any $m=1, \dots, g$, the following identities hold
\label{lemma_residues_g}
\begingroup
\allowdisplaybreaks
\begin{eqnarray}
\label{res3}
&& \frac{1}{2}\sum_{\substack{a_j\in B\\a_j\neq u_m}} \frac{\Omega(\paj)\varphi(\paj)}{a_j-u_m} +\underset{P=\pum}{\rm res} \frac{\Omega(P) \varphi(P)}{(u-u_m)du}=0\,;
\\
\label{res4}
&& \frac{1}{2}\sum_{\substack{a_j\in B\\a_j\neq u_m}} \Omega(\paj)W(\paj, \pum) +\underset{P=\pum}{\rm res} \frac{\Omega(P) W(P, \pum)}{du}=0\,;
\\
\label{res5}
&&  \Omega(P_0)v_m(P_0)+\Omega(P_1)v_m(P_1) +\sum_{i=1}^g\Omega(\pxi)v_m(\pxi) + \Omega(\pum)=0;
\\
\label{res6}
&&  \Omega(P_1)v_m(P_1) +\sum_{i=1}^gx_i\Omega(\pxi)v_m(\pxi) + u_m\Omega(\pum)=0.
\end{eqnarray}
\endgroup
\end{lemma}
{\it Proof.} All these relations are obtained by equating to zero the sums of residues of appropriate differentials on the compact surface $\Tx$. Namely,  \eqref{res3} corresponds to the sum of residues of the differential $\frac{\Omega(P) \varphi(P)}{(u-u_m)du}$. And \eqref{res4} corresponds to the differential $\frac{\Omega(P) W(P, \pum)}{du}$. For \eqref{res5} and \eqref{res6},  the differentials are respectively $\frac{\Omega(P)v_m(P)}{du(P)}$ and $\frac{u(P)\Omega(P)v_m(P)}{du(P)}$ where $v_m$ is the holomorphic differential \eqref{v} vanishing at every $\puj$ with $j\neq m.$ $\Box$

\section{Variation of the curve}
\label{sect_variation_pts}

\subsection{Rauch variation}
\label{sect_Rauch}

In this section we consider a family of hyperelliptic Riemann surfaces corresponding to the curves defined by equation \eqref{surf} where $x_1, \dots, x_g$ and $u_1, \dots, u_g$ are independent complex parameters varying in some small neighbourhood without coinciding. These variables parameterize complex structures of the curves through the standard local charts from Section \ref{sect_family}. Let us denote these surfaces by $\Txu.$ All the meromorphic differentials defined on our hyperelliptic Riemann surfaces depend  on the complex structure of the surface and thus vary under the variation of the branch points $x_j$ and $u_j.$ This variation is captured by the Rauch formulas \cite{Fay92}, see also \cite{KokoKoro}, which we give in this section. As before, we use the notation $a_j$ for an arbitrary branch point from the set $\B$ \eqref{B} although in this section we naturally assume that $a_j\notin\{0,1\}.$

A meromorphic differential $\Upsilon(P)$ defined on a surface from our family depends on the surface in question through the parameters $x_j$ and $u_j$ and on a point $P$ of that surface. The {\it Rauch variation} is the variation of the branch points $x_j$ and $u_j$ of the covering $u:\Tx\to \mathbb CP^1$  while keeping $u(P)$ fixed. Thus we define the {\it Rauch derivative} as follows for $a_k\in\B\setminus\{0,1\}$:
\begin{equation}
\label{Rauch-der}
\frac{\partial^{{\rm Rauch}}}{\partial a_k} \Upsilon(P) := \frac{\partial}{\partial a_k}\Big{|}_{u(P)=const} \Upsilon(P)\,.
\end{equation}

In the case of the Riemann bidifferential we need to require that both $u(P)$ and $u(Q)$ stay fixed. We have the following Rauch variational formula for the  $W$, see \cite{Fay92, KokoKoro}:
\begin{equation*}
\frac{\partial^{{\rm Rauch}}}{\partial a_k}W(P,Q) := \frac{\partial}{\partial a_k}\Big{|}_{\substack{u(P)=const\\ u(Q)=const}}W(P,Q) = \frac{1}{2}W(P, \pak)W(\pak, Q).
\end{equation*}
Assuming that the projections on the $u$-sphere of the chosen canonical homology basis do not change under the variation of branch points,   the Rauch variation of the Riemann bidifferential implies the following Rauch formulas for the holomorphic normalized differentials $\omega_j$ due to \eqref{bW} and for the Riemann matrix $\mathbb B$ of $\Txu$  \eqref{B}:
\begin{equation}
\label{RauchB}
\frac{\partial^{{\rm Rauch}} \omega_j(P)}{\partial a_k} = \frac{1}{2}\omega_j(\pak) W(P, \pak),
\qquad
\frac{\partial^{{\rm Rauch}} \mathbb B_{ij}}{\partial a_k} = \pi \i\,\omega_j(\pak)\omega_i(\pak), \qquad i,j=1, \dots, g\,.
\end{equation}
%
These Rauch variational formulas applied to definition \eqref{Omega}, \eqref{OW} of the differential of the third kind $\Omega$, yield
\begin{equation}
\label{Rauch-Omega}
\frac{\partial^{{\rm Rauch}} \Omega(P)}{\partial a_k} = \frac{1}{2}\Omega(\pak) W(P, \pak)\,.
\end{equation}
Evaluate both sides of the equation at $P=\paj$ for some $a_j\in\B$ and $a_j\neq a_k,$ we obtain
\begin{equation}
\label{Rauch-Omegaaj}
\frac{\partial^{{\rm Rauch}} \Omega(\paj)}{\partial a_k} = \frac{1}{2}\Omega(\pak) W(\paj, \pak)
\end{equation}
and we may do similarly for other Rauch formulas.
It is important to note that Rauch variation does not allow us to express the derivative of $\Omega(\pak)$ with respect to $a_k$. This situation is not covered by our definition of the Rauch derivative as in \eqref{Rauch-der} the $u$-image of the argument of the differential has to stay fixed under the variation. One can also see  that the right hand side of  \eqref{Rauch-Omega} cannot be evaluated at $\pak$ as  the bidifferential $W$  has a pole when both arguments coincide. However, for our purposes, we need to express $\partial_{a_k}\Omega(\pak)$ as well, and we do so in Section \ref{sect_variation}.

\subsection{Variation of dependent branch points on the Toda families of curves}
\label{sect_ux}

Let us now consider a Toda family of the curves defined by equation \eqref{surf} and find an expression for derivatives of the dependent branch points $u_1, \dots, u_{g}$  with respect to the independent branch points $x_1, \dots, x_g$.
\begin{theorem}
\label{thm_umder}
 Let $T$ be a Toda family of Riemann surfaces $\Tx$ corresponding to the curves of the equation \eqref{surf} where  the branch points $u_1, \dots, u_{g}$ of the covering $u:\Tx\to \mathbb CP^1$  are functions of the independently varying branch points $x_1, \dots, x_g.$ Let the chosen canonical homology basis  $\{\a_1,\dots, \a_g; \b_1, \dots, \b_g\}$ be as in Definition \ref{def_Todafamily} and $\Omega$ and $v_m$ be differentials \eqref{Omega} and \eqref{v}. Then we have for the functions $u_m(x_1, \dots, x_g):$

\begin{equation}
\label{umder_g}
\frac{\partial u_m}{\partial x_i}=-\frac{ \Omega(\pxi)}{\Omega(\pum)} v_m(\pxi)\,.
\end{equation}
\end{theorem}
{\it Proof.} Let us differentiate \eqref{Abel} with respect to an independent branch point $x_i$ using the Rauch formulas from Section \ref{sect_Rauch} and using \eqref{OW} to denote the integrals of $W$. In this way we obtain
\begingroup
\allowdisplaybreaks
\begin{multline}
\label{forlater}
{\bf \omega}(\pxi)\Omega_{\infty^-\infty^+}(\pxi)
+\sum_{k=1}^{g}{\bf \omega}(\puk)\Omega_{\infty^-\infty^+}(\puk) \frac{\partial u_k}{\partial x_i}
 \\
 = 2\pi \i{\bf \omega}(\pxi){\bf \omega}^T(\pxi) \N + 2\pi \i\sum_{k=1}^{g}{\bf \omega}(\puk){\bf \omega}^T(\puk) \N\frac{\partial u_k}{\partial x_i}\,,
\end{multline}
\endgroup
or equivalently, using definition \eqref{Omega} for $\Omega,$
\begin{equation}
\label{tempdiff_g}
\sum_{k=1}^{g}{\bf \omega}(P_{u_k})\Omega(P_{u_k}) \frac{\partial u_k}{\partial x_i}=-{\bf \omega}(P_{x_i})\Omega(P_{x_i})\,.
\end{equation}
Both sides of this equation are  $g$-component column vectors because $\omega$ is a vector of holomorphic normalized differentials. Thus, this is a linear system of equations for the $g$ unknown functions $\frac{\partial u_m}{\partial x_i}$, the matrix of coefficients being the $g\times g$ matrix $[\omega(P_{u_1})\Omega(P_{u_1}), \dots, \omega(P_{u_{g}})\Omega(P_{u_g})]$. Solving this system by Cramer's rule, we obtain
\begin{equation}
\label{umder}
 \frac{\partial u_m}{\partial x_i}=-\frac{ \Omega(\pxi)}{\Omega(\pum)} \frac{{\rm det}\,\mathcal M_m}{{\rm det}\,\mathcal M},
\end{equation}
where $\mathcal M$ stands for the $g\times g$ matrix $[\omega(P_{u_1}), \dots, \omega(P_{u_{g}})]$  and $\mathcal M_m$ is the matrix obtained from $\mathcal M$ by replacing the $m$th column with $\omega(\pxi).$ Note that expressions \eqref{tempdiff_g} and \eqref{umder} are invariant under replacing the column vector $\omega$ by a column vector whose components form any basis in the space of holomorphic one-forms on the surface $\Tx$. Thus by replacing  the vector $\omega$ by the column vector $v=(v_1, \dots, v_g)^T$ of the differentials \eqref{v}, we have the identity matrix instead of the matrix $\mathcal M$ in \eqref{umder}  and thus \eqref{umder} turns into \eqref{umder_g}
$\Box$

\begin{remark}
Note that in particular, for a continuous family of supports of the Pell equation of the form $[0,1]\cup(\cup_{g=1}^g[x_j, u_j])$ parameterized by $x_1, \dots, x_g$,   formula \eqref{umder_g} expresses derivatives of the positions of right-end points of the intervals with respect to the positions of their left-end points.
\end{remark}

\subsection{Variational formulas for the Toda families of curves}
\label{sect_variation}

Here we again consider a family of Riemann surfaces $\Tx, \;\x\in X$ corresponding to the curves of the equation  \eqref{surf} forming a Toda family in the sense of Definition \ref{def_Todafamily}. In particular, the branch points $u_1, \dots, u_{g}$  of the coverings $u:\Tx\to \mathbb CP^1$  are functions of the independently varying branch points $x_1, \dots, x_{g}.$
Our goal is to obtain derivatives of various meromorphic differentials defined on the family  $\Tx, \;\x\in X$ with respect to $x_1, \dots, x_{g}.$

The first such formulas can be easily obtained for  $\varphi(P_{a_j})$ and $ v_m(P_{a_j})$, the holomorphic differentials $\varphi$ \eqref{phi} and $v_m$ \eqref{v} evaluated at ramification points $\paj$ according to \eqref{evaluation}.
Let us assume that $u_\alpha$ stands for a dependent branch poing and $a_j$ denotes any of the independently varying branch points different from $x_i$.
The following formulas are obtained by a straightforward differentiation of \eqref{phi-eval}
\begin{equation}
\label{der_phi_aj_g}
\frac{\partial \varphi(P_{a_j})}{\partial x_i} = \frac{\varphi(P_{a_j})}{2}\left( \frac{1}{a_j-x_i} + \sum_{\alpha=1}^{g}\frac{1}{a_j-u_\alpha} \frac{\partial u_\alpha}{\partial x_i} \right),
\end{equation}
\begin{equation}
\label{der_phi_um}
\frac{\partial \varphi(P_{u_m})}{\partial x_i} = \frac{\varphi(P_{u_m})}{2}\left( \frac{1}{u_m-x_i} -\frac{\partial u_m}{\partial x_i}\sum_{\substack{a_j\in B\\a_j\neq u_m}} \frac{1}{u_m-a_j}  +\sum_{\substack{\alpha=1\\ \alpha \neq m}}^{g} \frac{1}{u_m-u_\alpha} \frac{\partial u_\alpha}{\partial x_i}  \right).
\end{equation}
Similarly, the evaluation of $v_m$ \eqref{v} at $\paj$ with any $a_j\in B$ gives
\begin{equation}
\label{vmpaj}
v_m(\paj)=\frac{\varphi (\paj)\prod_{\beta=1, \beta\neq m}^{g} (a_j-u_\beta)}{\varphi (P_{u_m})\prod_{\alpha=1, \alpha\neq m}^{g} (u_m-u_\alpha)}, \qquad   m=1, \cdots , g\,,
\end{equation}
which can be differentiated in a straightforward way using the obtained derivatives for $\varphi(\paj)$. This yields for $a_j\notin\{u_1, \dots, u_g, x_i\}$
\begingroup
\allowdisplaybreaks
\begin{multline}
\label{der_v_aj}
\frac{1}{ v_m(P_{a_j})}
\frac{\partial v_m(P_{a_j})}{\partial x_i} = \frac{1}{2}\left( \frac{1}{a_j-x_i} - \frac{1}{u_m-x_i} \right)+ \frac{1}{2}\frac{1}{a_j-u_m}\frac{\partial u_m}{\partial x_i} - \frac{1}{2} \sum_{\substack{\alpha=1\\ \alpha\neq m}}^{g}\frac{1}{a_j-u_\alpha} \frac{\partial u_\alpha}{\partial x_i}
 \\
+\frac{1}{2}\frac{\partial u_m}{\partial x_i} \left( \frac{1}{u_m} + \frac{1}{u_m-1} + \sum_{\alpha=1}^g\frac{1}{u_m-x_\alpha} \right) - \frac{1}{2} \sum_{\substack{\alpha=1\\ \alpha \neq m}}^g \frac{1}{u_m-u_\alpha} \left( \frac{\partial u_m}{\partial x_i} - \frac{\partial u_\alpha}{\partial x_i} \right)
\end{multline}
\endgroup
and
\begin{multline}
\label{der_v_xi}
\frac{1}{ v_m(\pxi)}\frac{\partial v_m(\pxi)}{\partial x_i}
 = - \frac{1}{2}\left( \frac{1}{x_i}+\frac{1}{x_i-1}+\sum_{\substack{\alpha=1\\ \alpha\neq i}}^g\frac{1}{x_i-x_\alpha}\right)  + \frac{1}{2}\sum_{\substack{\alpha=1\\ \alpha\neq m}}^g \frac{1}{x_i-u_\alpha}\left( 1 - \frac{\partial u_\alpha}{\partial x_i}\right)
\\
+ \frac{1}{2} \left(  \frac{1}{u_m} + \frac{1}{u_m-1} + \sum_{\substack{\alpha=1\\\alpha\neq i}}^g \frac{1}{u_m-x_\alpha}\right)\frac{\partial u_m}{\partial x_i}
- \frac{1}{2}\sum_{\substack{\alpha=1\\ \alpha \neq m}}^{g} \frac{1}{u_m-u_\alpha} \left(\frac{\partial u_m}{\partial x_i}-\frac{\partial u_\alpha}{\partial x_i}\right)  .
\end{multline}
In order to obtain further variational formulas and for the rest of the paper, we need the simple rational identities proven in the next lemma.
\begin{lemma}
\label{lemma_rational}
Let $M\in\mathbb N$ and a set $\{u_\alpha\}_{\alpha=1}^{M}\subset \mathbb C$ be given. For and any $x$  not in the set $\{u_\alpha\}_{\alpha=1}^{M}$, we have
\begin{equation}
\label{rat1}
\sum_{k=1}^{M} \frac{1}{(x-u_k)\prod_{ \alpha=1, \alpha \neq k}^M(u_k-u_\alpha)}=\frac{1}{\prod_{\alpha =1}^{M}(x-u_{\alpha})};
\end{equation}
\begin{multline}
\label{rat2}
\sum_{\substack{k=1\\ k\neq m}}^{M}\frac{1}{(u_k-x)(u_k-u_m)\prod_{ \alpha=1, \alpha \neq k}^M(u_k-u_\alpha)}
=
\frac{\sum_{\substack{\alpha=1\\ \alpha\neq m}}^M \frac{1}{u_m-u_\alpha}}{(u_m-x)\prod_{k=1, k\neq m}^M(u_m-u_k)}
\\
+ \frac{1}{(u_m-x)\prod_{k=1}^{M}(x-u_k)}+ \frac{1}{(u_m-x)^2\prod_{k=1, k\neq m}^M(u_m-u_k)}.
\end{multline}
\end{lemma}
{\it Proof.} Identity \eqref{rat1} is proven by comparing the poles and residues of the right and left hand sides. To establish identity \eqref{rat2}, let us first rewrite the left hand side as follows:
\begingroup
\allowdisplaybreaks
\begin{multline*}
\frac{1}{u_m-x}\sum_{\substack{k=1\\ k\neq m}}^{M}\left(\frac{1}{u_k-u_m}-\frac{1}{u_k-x}\right)\frac{1}{\prod_{\alpha\neq k}(u_k-u_\alpha)}
=\frac{1}{u_m-x}\frac{\partial}{\partial u_m}\sum_{\substack{k=1\\ k\neq m}}^{M}\frac{1}{(u_k-u_m)\prod_{\alpha\neq k,m}(u_k-u_\alpha)}
\\
+
\frac{1}{u_m-x}\sum_{k=1}^{M}\frac{1}{(x-u_k)\prod_{\alpha\neq k}(u_k-u_\alpha)}
+\frac{1}{(u_m-x)^2\prod_{\alpha\neq m}(u_m-u_\alpha)}.
\end{multline*}
\endgroup
It remains to use \eqref{rat1} in the first sum relative to the set $\{u_\alpha\}_{\substack{\alpha=1\\ \alpha \neq m}}^{M}\subset \mathbb C$ and $x=u_m$ as well as in the second sum.
$\Box$

Now, let us differentiate the quantity $\Omega(\paj)$ with $a_j$ being an independently varying branch point of the coverings $u:\Tx\to \mathbb CP^1.$

\begin{proposition}
\label{proposition_Omega_aj}
Let $\T\to  X$ be a Toda family of Riemann surfaces $\Tx$ as in Definition \ref{def_Todafamily} corresponding to the curves of equation \eqref{surf} where $\x=(x_1, \dots, x_g)$ are the independently varying branch points of the coverings $(\Tx,u)$ and the set of all branch points is $\B$ \eqref{branch}. Let $\Omega$ and $\varphi$ be, respectively, the meromorphic differential of the third kind defined by \eqref{Omega} and the holomorphic differential defined by \eqref{phi} on $\Tx$ for each $\x\in X$. Then for any $a_j\in \B\setminus\{x_i, u_1, \dots, u_g\}$ and $i=1, \dots, g$, we have for the evaluation of $\Omega$ at a ramification point $\paj$ obtained with respect to the standard local parameters $\sqrt{u(P)-a_j}$ according to \eqref{evaluation}:
\begin{equation*}
\frac{\partial \Omega(\paj)}{\partial x_i }  =\frac{\Omega(\pxi)\phi(\pxi)}{2(x_i-a_j)\phi(\paj)} \prod_{\alpha=1}^{g}\frac{x_i-u_\alpha}{a_j-u_\alpha} \,.
\end{equation*}
\end{proposition}
{\it Proof.}
From the Rauch formulas, we have
\begin{equation*}
\frac{\partial \Omega(\paj)}{\partial x_i }= \frac{1}{2} W(\pxi, \paj) \Omega(\pxi) + \frac{1}{2} \sum_{\alpha=1}^{g} W(\pualpha, \paj) \Omega(\pualpha) \frac{\partial u_\alpha}{\partial x_i}.
\end{equation*}
Now, let us use the following form for $W(P, \paj)$
\begin{equation}
\label{Waj}
W(P, \paj) = \frac{1}{u-a_j} \frac{\phi(P)}{\phi(\paj)} - \sum_{k=1}^g\beta_k^{(a_j)} \omega_k(P)
\end{equation}
where $\omega_k$ are the holomorphic normalized differentials and $\beta_k^{(a_j)}$ are normalization constants ensuring that the $a$-periods of $W(P, \paj)$ vanish. Plugging this into our derivative and using \eqref{Waj}  with $P=\pxi$ and $P=\pualpha$,
we see that the terms containing the normalization constants $\beta_k^{(a_j)}$  cancel out due to  \eqref{tempdiff_g}. We obtain
\begin{equation*}
\frac{\partial \Omega(\paj)}{\partial x_i }= \frac{1}{2} \frac{ \Omega(\pxi)}{x_i-a_j} \frac{\phi(\pxi)}{\phi(\paj)}
+ \frac{1}{2} \sum_{\alpha=1}^{g} \frac{\Omega(\pualpha)}{u_\alpha-a_j} \frac{\phi(\pualpha)}{\phi(\paj)}  \frac{\partial u_\alpha}{\partial x_i} .
\end{equation*}
Using \eqref{umder_g} for the derivatives of $u_\alpha$, we have
\begin{equation*}
\frac{\partial \Omega(\paj)}{\partial x_i }=  \frac{ \Omega(\pxi)\varphi(\pxi)}{2\varphi(P_{a_j})}
\left(\frac{1}{x_i-a_j}
- \sum_{\alpha=1}^{g} \frac{1}{u_\alpha-a_j}   \prod_{\beta=1, \beta\neq \alpha}^g \frac{x_i-u_\beta}{u_\alpha-u_\beta} \right).
\end{equation*}
Rewriting the sum over $\alpha$ slightly to facilitate the next step,
\begin{equation*}
\frac{\partial \Omega(\paj)}{\partial x_i }=  \frac{ \Omega(\pxi)\varphi(\pxi)}{2(x_i-a_j)\varphi(P_{a_j})}
\left(1
- \sum_{\alpha=1}^{g} \left(\frac{1}{x_i-u_\alpha} - \frac{1}{a_j-u_\alpha} \right)   \frac{ \prod_{s=1}^g(x_i-u_s)}{\prod_{\beta=1, \beta\neq \alpha}^g(u_\alpha-u_\beta)}\right)
\end{equation*}
and then applying identity \eqref{rat1}  with $x=x_i$ and $x=a_j$, we finish the proof of the proposition.
$\Box$

The next step is to differentiate the quantity obtained by evaluation of
the differential $\Omega$ at a ramification point whose position is not an independent variable. In this case, when applying the chain rule, we  need the following lemma, to obtain a derivative of $\Omega(\pum)$ with respect to $u_m$.
\begin{lemma}
\label{lemma_epsilon_g}
For $\mathcal B$ as in \eqref{branch},  let $\T_{\mathcal B}$ be a  family of Riemann surfaces  corresponding to the curves of equation \eqref{surf},  where we assume that all variable branch points  of the coverings $(\T_{\mathcal B},u)$ vary independently of each other.  Let $\Omega$ be the meromorphic differential of the third kind defined by \eqref{Omega} and $W$ be the fundamental Riemann bidifferential defined in Section \ref{sect_differentials}. Then for any $a_k\in \B$ , the evaluation of $\Omega$ at a ramification point $\pak$ obtained with respect to the standard local parameters $\sqrt{u(P)-a_k}$ according to \eqref{evaluation} satisfies

\begin{equation*}
\frac{\partial \Omega(P_{a_k})}{\partial a_k} = - \frac{1}{2}\sum_{\substack{a_j\in B \\ a_j\neq a_k}}\Omega(\paj)W(\pak,\paj) \,.
\end{equation*}
\end{lemma}
{\it Proof.} 
Let us consider the family of Riemann surfaces $\T_{\hat{\mathcal B}}$ corresponding to hyperelliptic curves of the equation
\begin{equation}
\label{hyperaj}
v^2=\prod_{j=1}^{2g+2}(u-a_j),
\end{equation}
parameterized by the elements of the set $\hat{\mathcal B}:=\{a_1, \dots, a_{2g+1}\}$. The branch points  $a_j$ of the covering $(\T_{\hat{\mathcal B}},u)$ vary continuously in some open set  $D\in\mathbb C^{2g+1}$, for which we assume that i) the points of $D$ are vectors with distinct components and ii) a canonical homology basis can be chosen for all the surfaces of the family $\T_{\hat{\mathcal B}}$ with $(a_1, \dots, a_{2g+1})\in D$ in such a way that their $u$-images are fixed for all $(a_1, \dots, a_{2g+1})\in D$, that is independent of $a_j$.

Let now $\varepsilon$ be a complex number with small absolute value such that for some fixed $(a_1^o, \dots, a^o_{2g+1})\in D$ we have $(a_1^o+\varepsilon, \dots, a^o_{2g+1}+\varepsilon)\in D$ and
 and consider  the covering $u_\varepsilon:\T_{\hat{\mathcal B}}\to\mathbb CP^1$  defined by $u_\varepsilon(P) = u(P)+\varepsilon$.
The covering $u_\varepsilon$  is ramified at the  ramification points $\paj=(a_j,0)$, $\;j=1, \dots, 2g+1,$ which are zeros of $du_\varepsilon,$ and at $P_\infty=(\infty, \infty)$ for any $\varepsilon$;  the corresponding branch points  are at $a_1+\varepsilon, \dots, a_{2g+1}+\varepsilon,\, \infty$. %

 We may evaluate the differential $\Omega$ defined on the surface corresponding to the equation \eqref{hyperaj} at a ramification point $\pak$ with respect to the standard local parameter at $\pak$ induced either by the covering $u$ or by the covering $u_\varepsilon.$ However, standard local parameters \eqref{coordinates}  are not affected by the shift $u\mapsto u+\varepsilon.$ Therefore, the quantity $\Omega(\pak)$ is independent of $\varepsilon.$ On the other hand, $\Omega(\pak)$ depends on the complex structure of the surface and thus is a function of branch points $\{a_j\}$ of the covering $u$. In the same way, $\Omega(\pak)$ is a function of branch points $\{a_j+\varepsilon\}$ of the covering $u_\varepsilon$. This allows us to write
\begin{equation*}
0=\frac{d}{d\varepsilon}{\Big |}_{\varepsilon=0}\Omega(\pak)=\sum_{j=1}^{2g+1} \frac{\partial}{\partial (a_j+\varepsilon)}\Omega(P_{a_k}){\Big |}_{\varepsilon=0}=\sum_{j=1}^{2g+1} \frac{\partial^{\rm Rauch}}{\partial a_j}\Omega(\pak)=\frac{\partial \Omega(\pak)}{\partial a_k}+\sum_{\substack{j=1\\ j\neq k}}^{2g+1} \frac{\partial^{\rm Rauch}}{\partial a_j}\Omega(\pak).
\end{equation*}
This, together with the Rauch formulas \eqref{Rauch-Omegaaj} proves the lemma for any set $\hat{\mathcal B}$ of branch points, and in particular for $\hat{\mathcal B}=\mathcal B$ \eqref{branch}.
$\Box$

Now we are ready for the next proposition which gives the derivative of $\Omega$ evaluated at a ramification point which is not included in Proposition \ref{proposition_Omega_aj}.
\begin{proposition}
\label{prop_Omega_um}
Let $ \T\to  X$ be a Toda family of Riemann surfaces $\Tx$ as in Definition \ref{def_Todafamily} corresponding to the curves of equation \eqref{surf} where $\x=(x_1, \dots, x_g)$ are the independently varying branch points of the coverings $u:\Tx\to\mathbb CP^1$ and the set of all branch points is $\B$ \eqref{branch}. Let $\Omega$ and $\varphi$ be the meromorphic differential of the third kind \eqref{Omega} and the holomorphic differential \eqref{phi} on each $\Tx$, respectively. Then for any $i, m=1, \dots, g$, we have for the evaluation of $\Omega$ at a ramification point $\pum$ obtained with respect to the standard local parameter $\sqrt{u(P)-u_m}$ according to \eqref{evaluation}:
\begin{equation*}
\frac{\partial \Omega(P_{u_m})}{\partial x_i} =\frac{\Omega(P_{u_m})}{2}
 \left( -\frac{1}{x_i-u_m} + \sum_{\substack{ \alpha=1\\ \alpha\neq m}}^{g}\frac{1}{u_m-u_\alpha}  -\frac{1}{\Omega(\pum)\varphi(\pum)}\sum_{\substack{a_j\in B\\a_j\neq u_m}} \frac{\Omega(P_{a_j})\varphi(P_{a_j})}{a_j-u_m}      \right)\frac{\partial u_m}{\partial x_i}\,.
\end{equation*}
\end{proposition}
{\it Proof.}
Differentiating $\Omega(P_{u_m})$ with respect to $x_i$ by the chain rule with the help of Rauch  formulas \eqref{Rauch-Omegaaj} and Lemma \ref{lemma_epsilon_g}, we have
\begin{multline*}
\frac{\partial \Omega(P_{u_m})}{\partial x_i} = \frac{\partial^{{\rm Rauch}} \Omega(P_{u_m})}{\partial x_i}
 + \sum_{\substack{k=1 \\ k\neq m}}^{g-1}\frac{\partial^{{\rm Rauch}} \Omega(P_{u_m})}{\partial u_k} \frac{\partial u_k}{\partial x_i} + \frac{\partial \Omega(P_{u_m})}{\partial u_m} \frac{\partial u_m}{\partial x_i}
 \\
 =\frac{1}{2}\Omega(\pxi)W(\pxi, \pum)  + \frac{1}{2}\sum_{\substack{k=1 \\ k\neq m}}^{g-1} \Omega(\puk) W(\puk, \pum) \frac{\partial u_k}{\partial x_i}
 -\frac{1}{2}\sum_{\substack{a_j\in B \\ a_j\neq u_m}} \Omega(\paj) W(\paj, \pum) \frac{\partial u_m}{\partial x_i}\,.
\end{multline*}
Let us now represent  $W(\cdot, \pum)$ similarly to \eqref{Waj} but with the help of the basis \eqref{v} in the space of holomorphic differentials as follows. For any $a_j\in B$ we have
\begin{equation}
\label{Wv}
W(P, \paj) = \frac{1}{u-a_j} \frac{\varphi(P)}{\varphi(\paj)} - \sum_{\alpha=1}^g\gamma_\alpha^{(a_j)} v_\alpha(P),
\end{equation}
where, again, $\gamma_\alpha^{(a_j)}$ are normalizing constants ensuring that the $a$-periods of  $W(\cdot, \paj)$ vanish.  Using this with $a_j=u_m,$ and evaluating $P$ at other ramification points $P=\pxi$, $P=\paj$, $P=\puk$ with $k\neq m$, using at the same time the defining conditions \eqref{vcond} for the differentials $v_\alpha$ as well as \eqref{res4}, the above expression for our derivative becomes
\begin{multline}
  \label{Omegau_temp}
\frac{\partial \Omega(P_{u_m})}{\partial x_i}
 =\frac{\Omega(\pxi)}{2}\left( \frac{1}{x_i-u_m} \frac{\varphi(\pxi)}{\varphi(\pum)} - \sum_{\alpha=1}^g\gamma_\alpha^{(u_m)} v_\alpha(\pxi)\right)
  \\
  + \sum_{\substack{k=1\\k\neq m}}^{g} \frac{\Omega(\puk)}{2} \left( \frac{1}{u_k-u_m} \frac{\varphi(\puk)}{\varphi(\pum)} - \gamma_k^{(u_m)} \right) \frac{\partial u_k}{\partial x_i}
 +\underset{P=\pum}{\rm res} \frac{\Omega(P) W(P, \pum)}{du} \frac{\partial u_m}{\partial x_i}\,.
\end{multline}

Now we note that, multiplying \eqref{Wv} by $\frac{\Omega(P)}{du}$, and using the characteristic properties \eqref{vcond} of the differentials $v_j$ and $du(P)=2\zeta_{u_m}(P)\zeta_{u_m}(P)$ for $P\sim \pum$ and $\zeta_{u_m}$ being the standard local coordinate from \eqref{coordinates},  one has
\begin{equation*}
\frac{1}{\varphi(P_{u_m})}\underset{P=P_{u_m}}{\rm res} \frac{\Omega(P) \varphi(P)}{(u-u_m)du} -
\underset{P=P_{u_m}}{\rm res} \frac{\Omega(P) W(P, P_{u_m})}{du} =
\underset{P=P_{u_m}}{\rm res} \frac{\Omega(P)}{du} \sum_{k=1}^g \gamma_k^{(u_m)} v_k(P)
=\frac{1}{2}\Omega(\pum)\gamma_m^{(u_m)}\,.
\end{equation*}
Together with \eqref{res3}, this allows us to express the residue appearing in \eqref{Omegau_temp}.
Plugging in this result as well as the derivatives \eqref{umder_g} of $u_m$ into \eqref{Omegau_temp}, we obtain an expression without terms containing normalization constants $\gamma_k^{(u_m)}$:
\begin{equation*}
\frac{\partial \Omega(P_{u_m})}{\partial x_i}
 =  \frac{\Omega(\pxi)\varphi(\pxi)}{2(x_i-u_m)\varphi(\pum)}
 - \frac{\Omega(\pxi)}{2\varphi(\pum)}  \sum_{\substack{k=1\\k\neq m}}^{g}  \frac{v_k(\pxi)\varphi(\puk)}{u_k-u_m}
 +\frac{ \Omega(\pxi)v_m(\pxi)}{2\Omega(\pum)\varphi(P_{u_m})}\sum_{a_j\neq u_m} \frac{\Omega(\paj)\varphi(\paj)}{a_j-u_m}   \,.
\end{equation*}
Let us now
plug in explicit expressions \eqref{v}  for the differentials $v_k$,  and  evaluate the sum over $k$ using  identity \eqref{rat2} with $x=x_i$ and $N=g.$ We obtain
\begin{multline*}
\frac{\partial \Omega(P_{u_m})}{\partial x_i}
 =  \frac{\Omega(\pxi)\varphi(\pxi)}{2(x_i-u_m)\varphi(\pum)} \left(
 -
\sum_{\beta\neq m} \frac{1}{u_m-u_\beta}
- \frac{1}{u_m-x_i}
\right.
\\
\left.
 +\frac{ 1}{\Omega(\pum)\varphi(\pum)} \sum_{a_j\neq u_m} \frac{\Omega(P_{a_j})\varphi(P_{a_j})}{a_j-u_m}   \right)
 \frac{\prod_{\alpha=1}^{g} (x_i-u_\alpha) }{\prod_{\alpha=1, \alpha\neq m}^g(u_m-u_\alpha)}\,,
\end{multline*}
which is equivalent to the statement of the proposition due to expression \eqref{umder_g} for $\partial_{x_i}u_m$.
$\Box$

\begin{proposition}
\label{prop_Omega_xk} In the situation of Proposition \ref{prop_Omega_um},
for any $k=1, \dots, g$, we have for the evaluation of $\Omega$ at a ramification point $\pxk$ obtained with respect to the standard local parameter $\sqrt{u(P)-x_k}$ according to \eqref{evaluation}:
\begin{equation*}
\frac{\partial \Omega(\pxk)}{\partial x_k} =-\frac{\Omega(\pxk)}{2}
 \left(  \sum_{\substack{ \alpha=1}}^{g}\frac{1}{u_\alpha-x_k} \prod_{\substack{\beta=1\\ \beta\neq \alpha}}^g\frac{x_k-u_\beta}{u_\alpha-u_\beta}  + \sum_{\substack{a_j\in B\\a_j\neq x_k}}\frac{1}{a_j-x_k} \frac{\Omega(P_{a_j})\varphi(P_{a_j})}{\Omega(\pxk)\varphi(\pxk)}     \right)\,.
\end{equation*}
\end{proposition}
{\it Proof.} Lemma \ref{lemma_epsilon_g} and the chain rule imply
\begin{equation}
\label{temppr3}
\frac{\partial \Omega(\pxk)}{\partial x_k} = - \frac{1}{2}\sum_{\substack{a_j\in B\\a_j\neq x_k}} \Omega(\paj)W(\paj, \pxk) + \sum_{\substack{j=1}}^g\frac{\partial^{\rm Rauch}\Omega(\pxk)}{\partial u_j}\frac{\partial u_j}{\partial x_k}.
\end{equation}
Now we replace the first sum by the residue at $\pxk$ of $\frac{\Omega(P)W(P, \pxk)}{du}$ analogously to relation \eqref{res4} and use the Rauch formulas \eqref{Rauch-Omega} in the second sum as well as \eqref{umder_g} for $\partial_{x_k}u_j$. In doing so, we write the bidifferential $W$ as in \eqref{Wv} replacing $a_j$ by $x_k$. This yields for the residue:
\begin{equation*}
\underset{P=\pxk}{\rm res}\frac{\Omega(P) W(P, \pxk)}{du}=\frac{1}{\phi(\pxk)}\underset{P=\pxk}{\rm res}\frac{\Omega(P) \varphi(P)}{(u-x_k)du} - \frac{1}{2}\sum_{j=1}^g \Omega(\pxk) v_j(\pxk)\gamma_j^{(x_k)}.
\end{equation*}
The sum over $j$ cancels against the terms coming out of the sum over $j$ in \eqref{temppr3}
\begin{equation*}
\sum_{\substack{j=1}}^g\frac{\partial^{\rm Rauch}\Omega(\pxk)}{\partial u_j}\frac{\partial u_j}{\partial x_k}
=-\frac{\Omega(\pxk)}{2}\sum_{j=1}^g v_j(\pxk) \left( \frac{1}{u_j-x_k} \frac{\varphi(\puj)}{\varphi(\pxk)} - \gamma_j^{(x_k)}\right).
\end{equation*}
Replacing now the remaining residue through \eqref{res3} which holds upon replacing $u_m$ by $x_k$, and replacing $v_j(\pxk)$ by their explicit representation \eqref{v},
we finish the proof.
$\Box$

\section{Isoperiodicity equations}
\label{sect_isoeqs}

In this section, we derive a system of differential equations for the positions of dependent branch points with respect to the independent ones of the coverings defining a Toda family of Riemann surfaces. In Section \ref{sect_generalized}, we interpret this system as isoperiodic deformations of a Riemann surface carrying a differential of the third kind with fixed periods on it. Recall that existence of Toda families follows from Theorems 2.7 and 2.12 from \cite{PS1999}; we also prove it in Section \ref{sect_generalized}.
\begin{theorem}
\label{thm_support}
Let $\T\to  X$ be a Toda family of hyperelliptic Riemann surfaces $\Tx$  corresponding to the curves of equation \eqref{surf} parameterized by $\x=(x_1, \dots, x_g)\in X$ as in Definition \ref{def_Todafamily}. Then the  branch points $(u_1, \dots, u_g)$ of the coverings $u:\Tx\to\mathbb CP^1$ satisfy the following system with respect to variables $\x=(x_1, \dots, x_g)$.

For $m,k,n\in\{1, \dots, g\}$ with $k\neq n$
\begingroup
\allowdisplaybreaks
\begin{multline}
\label{um_xkxn}
\frac{\partial^2 u_m}{\partial x_k\partial x_n}
=\frac{1}{2}\left( \frac{1}{x_k-x_n} - \frac{1}{u_m-x_n} \right)\frac{\partial u_m}{\partial x_k}
+\frac{1}{2} \left( \frac{1}{x_n-x_k}-\frac{1}{u_m-x_k}  \right)\frac{\partial u_m}{\partial x_n}
\\
+ \frac{1}{2}\left(\frac{2}{u_m} + \frac{2}{u_m-1} + \sum_{\substack{j=1\\j\neq k,m}}^g\frac{1}{u_m-x_j}-\sum_{\substack{j=1 \\ j\neq m}}^g \frac{1}{u_m-u_j}     \right)\frac{\partial u_m}{\partial x_n}\frac{\partial u_m}{\partial x_k}
\\
+ \frac{1}{4} \frac{\partial u_m}{\partial x_k}\sum_{\substack{\alpha=1\\ \alpha\neq m}}^{g}\left( \frac{1}{u_m-u_\alpha}-\frac{1}{x_k-u_\alpha}\right) \frac{\partial u_\alpha}{\partial x_n}
+ \frac{1}{4} \frac{\partial u_m}{\partial x_n}\sum_{\substack{\alpha=1\\ \alpha\neq m}}^{g}\left( \frac{1}{u_m-u_\alpha}-\frac{1}{x_n-u_\alpha}\right) \frac{\partial u_\alpha}{\partial x_k}
\\
-\frac{1}{2} \frac{\partial u_m}{\partial x_k} \frac{\partial u_m}{\partial x_n} \sum_{i=1}^g \frac{\partial u_m}{\partial x_i}  \left( \sum_{\substack{j=1 \\ j\neq m}}^g \frac{1}{u_m-u_j} +\frac{1}{u_m}   \right)
%
 -\frac{1}{2u_m(u_m-1)}  \frac{\partial u_m}{\partial x_k} \frac{\partial u_m}{\partial x_n}   \sum_{i=1}^g x_i\frac{\partial u_m}{\partial x_i}
\\
-  \frac{1}{2} \frac{\partial u_m}{\partial x_k} \frac{\partial u_m}{\partial x_n} \sum_{j=1}^g\frac{1}{\prod_{\alpha=1}^g(x_j-u_\alpha)} \frac{\partial u_m}{\partial x_j}\prod_{\substack{\beta=1\\ \beta\neq  m}}^g(u_m-u_\beta)
\\
-
\frac{1}{2} \frac{\partial u_m}{\partial x_k} \frac{\partial u_m}{\partial x_n} \sum_{\substack{j=1 \\ j\neq m}}^g
  \sum_{i=1}^g\frac{\partial u_m}{\partial x_i} \left(\frac{1}{u_m-u_j}-\frac{1}{x_i-u_j}\right)  \frac{1}{\prod_{\substack{\alpha=1 \\ \alpha\neq j}}^g(u_j-u_\alpha)}\prod_{\substack{\beta=1 \\ \beta\neq m}}^g (u_m-u_\beta);
\end{multline}
\endgroup
and for any $m,k\in\{1, \dots, g\}$
\begingroup
\allowdisplaybreaks
\begin{multline}
\label{um_xkxk}
\frac{\partial^2 u_m}{\partial x_k^2}=\frac{1}{2}\left( \frac{u_m}{x_k-1} - \frac{u_m-1}{x_k} - \frac{1}{x_k-u_m} \right)
  -\frac{1}{2} \sum_{\substack{j=1\\j\neq k}}^g\left(\frac{1}{x_k-u_m}-\frac{1}{x_k-x_j}\right)\frac{\partial u_m}{\partial x_j}
\\
- \frac{1}{2}\frac{\partial u_m}{\partial x_k}\left( \frac{2}{x_k}+\frac{2}{x_k-1}+\sum_{\substack{\alpha=1\\ \alpha\neq k}}^g\frac{1}{x_k-x_\alpha}-\sum_{\substack{\alpha=1\\ \alpha\neq m}}^g \frac{1}{x_k-u_\alpha} -\frac{1}{x_k-u_m}\right)
\\
+\frac{1}{2}\left( \frac{1}{x_k}- \frac{1}{x_k-1}\right)\sum_{\substack{i=1\\ i\neq k}}^gx_i\frac{\partial u_m}{\partial x_i}
+\frac{1}{2}\left( \frac{1}{x_k-u_m}-\frac{1}{x_k}\right)\sum_{\substack{i=1\\ i\neq k}}^g\frac{\partial u_m}{\partial x_i}
\\
+ \frac{1}{2}\left( \frac{\partial u_m}{\partial x_k}\right)^2 \left(  \frac{2}{u_m} + \frac{2}{u_m-1} + \sum_{\substack{\alpha=1\\\alpha\neq k}}^g \frac{1}{u_m-x_\alpha}-\sum_{\substack{\alpha=1\\ \alpha \neq m}}^{g} \frac{1}{u_m-u_\alpha} +\frac{1}{x_k-u_m} \right)
\\
+\frac{1}{2}\frac{\partial u_m}{\partial x_k}\sum_{\substack{\alpha=1\\ \alpha \neq m}}^{g} \left(\frac{1}{u_m-u_\alpha} -\frac{1}{x_k-u_\alpha} \right)\frac{\partial u_\alpha}{\partial x_k}
-\frac{1}{2}\left(\frac{\partial u_m}{\partial x_k}\right)^2\sum_{j=1}^g\frac{1}{\prod_{\alpha=1}^g(x_j-u_\alpha)} \frac{\partial u_m}{\partial x_j}\prod_{\substack{\beta=1\\ \beta\neq  m}}^g(u_m-u_\beta)
\\
-\frac{1}{2}\left(\frac{\partial u_m}{\partial x_k}\right)^2\sum_{i=1}^g \frac{\partial u_m}{\partial x_i}  \left(\frac{1}{u_m}+\sum_{\substack{j=1 \\ j\neq m}}^g \frac{1}{u_m-u_j}\right)
 -\frac{1}{2u_m(u_m-1)}\left(\frac{\partial u_m}{\partial x_k}\right)^2   \sum_{i=1}^gx_i\frac{\partial u_m}{\partial x_i}
\\
-
\frac{1}{2}\left(\frac{\partial u_m}{\partial x_k}\right)^2\sum_{\substack{j=1 \\ j\neq m}}^g
  \sum_{i=1}^g\frac{\partial u_m}{\partial x_i} \left(\frac{1}{u_m-u_j}-\frac{1}{x_i-u_j}\right)  \frac{1}{\prod_{\substack{\alpha=1 \\ \alpha\neq j}}^g(u_j-u_\alpha)}\prod_{\substack{\beta=1 \\ \beta\neq m}}^g (u_m-u_\beta).
\end{multline}
\endgroup

\end{theorem}
\begin{corollary}
\label{cor_Chebyshev}
Let $x_1<x_2< \dots< x_g$ be independent real variables and $(u_1, \dots, u_g)$ be functions of $x_1, \dots, x_g$ such that $x_j<u_j$ and $[0,1]\cup(\cup_{g=1}^g[x_j, u_j])\subset \mathbb R$ is a continuous family of supports of the Pell equation \eqref{Pell}. Then the functions $u_m(x_1, \dots, x_g)$ satisfy system \eqref{um_xkxn}, \eqref{um_xkxk}.
\end{corollary}
{\it Proof.} This follows from Theorem \ref{thm_support}, Definition \ref{def_Todafamily} of a Toda family, of which the situation of the corollary is a particular case, and Theorems 2.7 and 2.12 from \cite{PS1999} stating the existence of continuous families of supports of the Pell equation of the form in question. $\Box$

To prove Theorem \ref{thm_support}, we need an auxiliary result of the next lemma.
\begin{lemma}
\label{prop_Omega_pum}
For a Toda family  $\T\to  X$ of Riemann surfaces $\Tx$  corresponding to the curves of equation \eqref{surf}, let  $\Omega$ be the meromorphic differential of the third kind \eqref{Omega} and let $\varphi$ be the holomorphic differential \eqref{phi} on each surface from the family. Let the quantities $\phi(\pum)$ and $\Omega(\pum)$ be obtained according to \eqref{evaluation} with respect to the standard local parameter $\sqrt{u(P)-u_m}$ at the ramification point $\pum$ of the covering $u:\Tx\to\mathbb CP^1.$ Then for any $m=1, \dots, g$
\begingroup
\allowdisplaybreaks
\begin{multline}
\label{res_main}
\frac{1}{\Omega(\pum)\phi(\pum)}\underset{P=\pum}{\rm res} \frac{\Omega(P) \varphi(P)}{(u-u_m)du}
=\frac{1}{2}\sum_{j=1}^g\frac{1}{\prod_{\alpha=1}^g(x_j-u_\alpha)} \frac{\partial u_m}{\partial x_j}\prod_{\substack{\beta=1\\ \beta\neq  m}}^g(u_m-u_\beta)
\\
+\frac{1}{2}\left(\sum_{i=1}^g \frac{\partial u_m}{\partial x_i}  - 1 \right) \sum_{\substack{j=1 \\ j\neq m}}^g \frac{1}{u_m-u_j}
+\frac{1}{2u_m}\left(\sum_{i=1}^g \frac{\partial u_m}{\partial x_i}  - 1 \right)
 +\frac{1}{2(u_m-1)} \left(  \sum_{i=1}^g\frac{x_i}{u_m}\frac{\partial u_m}{\partial x_i} - 1    \right)
\\
+
\frac{1}{2}\sum_{\substack{j=1 \\ j\neq m}}^g
  \sum_{i=1}^g\frac{\partial u_m}{\partial x_i} \left(\frac{1}{u_m-u_j}-\frac{1}{x_i-u_j}\right)  \frac{1}{\prod_{\substack{\alpha=1 \\ \alpha\neq j}}^g(u_j-u_\alpha)}\prod_{\substack{\beta=1 \\ \beta\neq m}}^g (u_m-u_\beta).
\end{multline}
\endgroup
\end{lemma}
{\it Proof.} From \eqref{res3} of Lemma \ref{lemma_residues_g} we have
\begin{equation*}
\underset{P=\pum}{\rm res} \frac{\Omega(P) \varphi(P)}{(u-u_m)du}=\frac{1}{2}\frac{\Omega(P_0)\phi(P_0)}{u_m} -\frac{1}{2}\frac{\Omega(P_1)\phi(P_1)}{1-u_m}
-\frac{1}{2}\sum_{j=1}^g\frac{\Omega(\pxj)\phi(\pxj)}{x_j-u_m}
-\frac{1}{2}\sum_{\substack{j=1 \\ j\neq m}}^g\frac{\Omega(\puj)\phi(\puj)}{u_j-u_m}.
\end{equation*}
Next, we divide by ${\Omega(\pum)\phi(\pum)}$ and substitute the expression for $\Omega(\puj)$  obtained from \eqref{res5} with $m=j$. We use \eqref{umder_g} to rewrite parts of the expression in terms of derivatives of $u_m.$ This yields
\begin{multline}
\label{temp_prop}
\frac{1}{\Omega(\pum)\phi(\pum)}\underset{P=\pum}{\rm res} \frac{\Omega(P) \varphi(P)}{(u-u_m)du}=
\frac{1}{2}\frac{\Omega(P_0)\phi(P_0)}{\Omega(\pum)\phi(\pum)}\left( \frac{1}{u_m}-\sum_{\substack{j=1 \\ j\neq m}}^g\frac{ \prod_{\alpha=1}^g(-u_\alpha)  }{u_j(u_j-u_m)}\prod_{\substack{\beta=1 \\ \beta\neq j}}^g \frac{1}{(u_j-u_\beta)}\right)
\\
 +\frac{1}{2}\frac{\Omega(P_1)\phi(P_1)}{\Omega(\pum)\phi(\pum)}\left(\frac{1}{u_m-1} -\sum_{\substack{j=1 \\ j\neq m}}^g\frac{ \prod_{\alpha=1}^g(1-u_\alpha)  }{(u_j-1)(u_j-u_m)}\prod_{\substack{\beta=1 \\ \beta\neq j}}^g \frac{1}{(u_j-u_\beta)}\right)
\\
+\frac{1}{2}\sum_{j=1}^g\frac{1}{\prod_{\alpha=1}^g(x_j-u_\alpha)} \frac{\partial u_m}{\partial x_j}\prod_{\substack{\beta=1\\ \beta\neq  m}}^g(u_m-u_\beta)
\\
-
\frac{1}{2}\sum_{\substack{j=1 \\ j\neq m}}^g
  \sum_{i=1}^g\frac{\partial u_m}{\partial x_i} \left(\frac{1}{u_j-u_m}+\frac{1}{x_i-u_j}\right)  \frac{1}{\prod_{\substack{\alpha=1 \\ \alpha\neq j}}^g(u_j-u_\alpha)}\prod_{\substack{\beta=1 \\ \beta\neq m}}^g (u_m-u_\beta).
\end{multline}
Now we apply rational identity \eqref{rat2} with $N=g$ for the sums in the first and second lines: for the first line we set $x=0$ and for the second line $x=1.$ This gives
\begin{multline*}
\frac{1}{\Omega(\pum)\phi(\pum)}\underset{P=\pum}{\rm res} \frac{\Omega(P) \varphi(P)}{(u-u_m)du}
=\frac{1}{2}\sum_{j=1}^g\frac{1}{\prod_{\alpha=1}^g(x_j-u_\alpha)} \frac{\partial u_m}{\partial x_j}\prod_{\substack{\beta=1\\ \beta\neq  m}}^g(u_m-u_\beta)
\\
-\frac{1}{2}\frac{\Omega(P_0)\phi(P_0)}{\Omega(\pum)\phi(\pum)}\left( \sum_{\substack{\alpha=1 \\ \alpha\neq m}}^g \frac{1}{u_m-u_\alpha}
+\frac{1}{u_m}
\right)\frac{\prod_{\alpha=1}^g(-u_\alpha) }{u_m\prod_{k\neq m}(u_m-u_k)}
\\
 -\frac{1}{2}\frac{\Omega(P_1)\phi(P_1)}{\Omega(\pum)\phi(\pum)}\left(    \sum_{\substack{\alpha=1 \\ \alpha\neq m}}^g \frac{1}{u_m-u_\alpha}
+ \frac{1}{u_m-1}  \right)\frac{\prod_{\alpha=1}^g(1-u_\alpha)}{(u_m-1)\prod_{k\neq m}(u_m-u_k)}
\\
-
\frac{1}{2}\sum_{\substack{j=1 \\ j\neq m}}^g
  \sum_{i=1}^g\frac{\partial u_m}{\partial x_i} \left(\frac{1}{u_j-u_m}+\frac{1}{x_i-u_j}\right) \frac{1}{\prod_{\substack{\alpha=1 \\ \alpha\neq j}}^g(u_j-u_\alpha)}
  \prod_{\substack{\beta=1 \\ \beta\neq m}}^g (u_m-u_\beta).
\end{multline*}

Let us now express $\frac{\Omega(P_1)\phi(P_1)}{\Omega(\pum)\phi(\pum)}$ from \eqref{res6} of Lemma \ref{lemma_residues_g} as follows
\begin{equation}
\label{1}
\frac{\Omega(P_1)\phi(P_1)}{\Omega(\pum)\phi(\pum)}= \left(  \sum_{i=1}^gx_i\frac{\partial u_m}{\partial x_i} - u_m    \right)\prod_{\substack{\alpha=1 \\ \alpha\neq m}}^g\frac{u_m-u_\alpha}{1-u_\alpha}.
\end{equation}
This  and \eqref{res5} of Lemma \ref{lemma_residues_g} allow to  express $\frac{\Omega(P_0)\phi(P_0)}{\Omega(\pum)\phi(\pum)}$ as follows
\begin{equation}
\label{0}
\frac{\Omega(P_0)\phi(P_0)}{\Omega(\pum)\phi(\pum)}= \left(\sum_{i=1}^g(1-x_i) \frac{\partial u_m}{\partial x_i} + u_m  - 1 \right)\prod_{\substack{\alpha=1 \\ \alpha\neq m}}^g \frac{u_\alpha-u_m}{u_\alpha}.
\end{equation}
Plugging \eqref{1} and \eqref{0} into \eqref{temp_prop}, we prove the lemma.
$\Box$

{\it Proof of Theorem \ref{thm_support}.}
Let us differentiate $\partial{u_m}/{\partial x_k}$ given by \eqref{umder_g} with respect to $x_n$ for $n\neq k.$ This gives

\begin{equation}
\label{der1-initial}
\frac{\partial^2 u_m}{\partial x_k\partial x_n}=\frac{\partial}{\partial x_n} \left\{  -  \frac{ v_{m}(\pxk)\Omega(\pxk)}{ \Omega(\pum)}   \right\}
=\frac{\partial u_m}{\partial x_k} \left( \frac{\frac{\partial v_m(\pxk)}{\partial x_n}}{v_m(\pxk)} + \frac{\frac{\partial \Omega(\pxk)}{\partial x_n}}{\Omega(\pxk)} -  \frac{\frac{\partial \Omega(\pum)}{\partial x_n}}{\Omega(\pum)}    \right).
\end{equation}

Proposition \eqref{proposition_Omega_aj} implies the following relation setting $a_j=x_k:$
\begin{equation}
\label{der-rel}
\frac{\partial u_m}{\partial x_k}\frac{1}{\Omega(\pxk)} \frac{\partial \Omega(\pxk)}{\partial x_n } = \frac{1}{2}\frac{\partial u_m}{\partial x_n} \left( \frac{1}{x_k-u_m} + \frac{1}{x_n-x_k} \right),
\end{equation}
which gives us the second term in the right hand side of \eqref{der1-initial}. The first  term is computed in \eqref{der_v_aj}.
For the third term, let us rewrite the expression obtained in Proposition \ref{prop_Omega_um}  using relation \eqref{res3}, that is
\begin{equation}
\label{OmegaPum_xn}
\frac{1}{\Omega(P_{u_m})}\frac{\partial \Omega(P_{u_m})}{\partial x_n} =\frac{1}{2}
 \left(  \sum_{\substack{ \alpha=1\\ \alpha\neq m}}^{g}\frac{1}{u_m-u_\alpha}-\frac{1}{x_n-u_m}   +\frac{2}{\Omega(\pum)\phi(\pum)}\underset{P=\pum}{\rm res} \frac{\Omega(P) \varphi(P)}{(u-u_m)du}      \right)\frac{\partial u_m}{\partial x_n}\,.
\end{equation}
A rational expression for the remaining residue is provided by \eqref{res_main} in Lemma \eqref{prop_Omega_pum}.
Plugging all these results into \eqref{der1-initial}, we obtain an equation that coincides with \eqref{um_xkxn} except for the third line in \eqref{um_xkxn}.
Instead of the third line in \eqref{um_xkxn}, we obtain the following non-symmetric expression
\begin{equation}
\label{asymm}
S:=\frac{1}{2} \frac{\partial u_m}{\partial x_k}\sum_{\substack{\alpha=1\\ \alpha\neq m}}^{g}\left( \frac{1}{u_m-u_\alpha}-\frac{1}{x_k-u_\alpha}\right) \frac{\partial u_\alpha}{\partial x_n}.
\end{equation}
Let us show that \eqref{asymm} is equivalent to the expression in the third line of \eqref{um_xkxn}. To this end, we rewrite the derivatives in \eqref{asymm} in terms of differentials $\Omega$ and $v_j$ using \eqref{umder_g}:
\begin{multline*}
S
=
\frac{1}{2} \frac{ \Omega(\pxk) v_m(\pxk)}{\Omega(\pum)}\sum_{\substack{\alpha=1\\ \alpha\neq m}}^{g} \frac{x_k-u_m}{(u_m-u_\alpha)(x_k-u_\alpha)} \frac{ \Omega(\pxn)v_\alpha(\pxn)}{\Omega(\pualpha)}
\\
=\frac{1}{2} \frac{ \Omega(\pxk) \Omega(\pxn)v_m(\pxk)v_m(\pxn)}{\Omega(\pum)}\sum_{\substack{\alpha=1\\ \alpha\neq m}}^{g} \frac{(x_k-u_m)(x_n-u_m)}{(u_m-u_\alpha)(x_k-u_\alpha)(x_n-u_\alpha)} \frac{ \phi(\pum)}{\Omega(\pualpha)\phi(\pualpha)} \frac{\prod_{\substack{j=1\\ j\neq m}}(u_m-u_j)}{\prod_{\substack{i=1\\ i\neq \alpha}}(u_\alpha-u_i)},
\end{multline*}
which proves that expression \eqref{asymm} is symmetric in $n$ and $k$ and thus can be replaced by its symmetrization in the third line in \eqref{um_xkxn}.

Let us now prove the second equation \eqref{um_xkxk}. We need to show that the following expression is rational in $u_j, \,x_j$ and $\partial_{x_i} u_j:$
\begin{equation}
\label{der2-initial}
\frac{\partial^2 u_m}{\partial x_k^2}=\frac{\partial}{\partial x_k} \left\{  -  \frac{ v_{m}(\pxk)\Omega(\pxk)}{ \Omega(\pum)}   \right\}
=\frac{\partial u_m}{\partial x_k} \left( \frac{\frac{\partial v_m(\pxk)}{\partial x_k}}{v_m(\pxk)} + \frac{\frac{\partial \Omega(\pxk)}{\partial x_k}}{\Omega(\pxk)} -  \frac{\frac{\partial \Omega(\pum)}{\partial x_k}}{\Omega(\pum)}    \right).
\end{equation}
Such a rational expression  is given by \eqref{der_v_xi} for the first term and, upon replacing $n$ by $k$, by the combination of \eqref{OmegaPum_xn}  and \eqref{res_main} for the last term. Thus it remains to evaluate the middle term in \eqref{der2-initial}. Let us take the expression for $\partial_{x_k}\Omega(\pxk)$ from Proposition \ref{prop_Omega_xk} and replace $\Omega(\puj)$ by their expression obtained from \eqref{res5}, namely, $\Omega(\puj)=-\Omega(P_0)v_j(P_0)-\Omega(P_1)v_j(P_1) -\sum_{i=1}^g\Omega(\pxi)v_j(\pxi)$. This yields
\begin{multline*}
\frac{1}{\Omega(\pxk)}\frac{\partial \Omega(\pxk)}{\partial x_k} =
\frac{1}{2} \frac{\Omega(P_0)\varphi(P_0)}{x_k\Omega(\pxk)\varphi(\pxk)} \left( 1- \left( \prod_{\beta=1}^g(-u_\beta)\right)\sum_{j=1}^g\left(\frac{1}{u_j-x_k}-\frac{1}{u_j}\right) \prod_{\substack{\alpha=1\\\alpha\neq j}}^g\frac{1}{u_j-u_\alpha} \right)
\\   +\frac{1}{2} \frac{\Omega(P_1)\varphi(P_1)}{(x_k-1)\Omega(\pxk)\varphi(\pxk)}    \left( 1 -\left(\prod_{\beta=1}^g(1-u_\beta)\right)\sum_{j=1}^g\left(\frac{1}{u_j-x_k}-\frac{1}{u_j-1}\right) \prod_{\substack{\alpha=1\\\alpha\neq j}}^g\frac{1-u_\alpha}{u_j-u_\alpha}\right)
\\
   -\frac{1}{2} \sum_{\substack{j=1\\j\neq k}}^g\frac{1}{x_j-x_k}\frac{\Omega(P_{x_j})\varphi(P_{x_j})}{\Omega(\pxk)\varphi(\pxk)}  \left( 1    -\left( \prod_{\beta=1}^g(x_j-u_\beta)\right) \sum_{i=1}^g\left(\frac{1}{u_i-x_k}    +\frac{1}{x_j-u_i}  \right)\prod_{\substack{\alpha=1\\\alpha\neq i}}^g\frac{1}{u_i-u_\alpha}\right)\,.
\end{multline*}
Applying now rational identity \eqref{rat1} in each line, we have
\begin{multline}
\label{Omegapxk_temp}
\frac{1}{\Omega(\pxk)}\frac{\partial \Omega(\pxk)}{\partial x_k} =
\frac{1}{2} \frac{\Omega(P_0)\varphi(P_0)}{x_k\Omega(\pxk)\varphi(\pxk)} \prod_{\alpha=1}^g\frac{(-u_\alpha)}{x_k-u_\alpha}
   +\frac{1}{2} \frac{\Omega(P_1)\varphi(P_1)}{(x_k-1)\Omega(\pxk)\varphi(\pxk)}   \prod_{\alpha=1}^g\frac{1-u_\alpha}{ x_k-u_\alpha}
\\
   -\frac{1}{2} \sum_{\substack{j=1\\j\neq k}}^g\frac{1}{x_j-x_k}\frac{\Omega(P_{x_j})\varphi(P_{x_j})}{\Omega(\pxk)\varphi(\pxk)} \prod_{\alpha=1}^g  \frac{x_j-u_\alpha}{x_k-u_\alpha}\,.
\end{multline}
Multiplying now by $\partial_{x_k}u_m$ on both sides and writing this derivative as in \eqref{umder_g} in the right hand side, we see that $\Omega(\pxk)\varphi(\pxk)$ in the denominator is replaced by $\Omega(\pum)$, which allows to write the terms in the last line of \eqref{Omegapxk_temp} using the derivatives $\partial_{x_i} u_m$. It remains   to replace the quantities $\frac{\Omega(P_0)\varphi(P_0)}{\Omega(\pum)}$ and $\frac{\Omega(P_1)\varphi(P_1)}{\Omega(\pum)}$ using \eqref{res5} and \eqref{res6} similarly to \eqref{0} and \eqref{1} to obtain the following rational expression for the middle term of \eqref{der2-initial}:
\begin{multline}
\label{middle}
\frac{1}{\Omega(\pxk)}\frac{\partial \Omega(\pxk)}{\partial x_k} \frac{\partial u_m}{\partial x_k} =
\frac{1}{2}\left( \frac{1}{x_k-u_m}-\frac{1}{x_k}\right)\left(\sum_{i=1}^g(1-x_i)\frac{\partial u_m}{\partial x_i} + u_m-1 \right)
\\   + \frac{1}{2} \left(\frac{1}{x_k-u_m} - \frac{1}{x_k-1} \right) \left(  \sum_{i=1}^gx_i\frac{\partial u_m}{\partial x_i} - u_m  \right)
  -\frac{1}{2} \sum_{\substack{j=1\\j\neq k}}^g\left(\frac{1}{x_k-u_m}-\frac{1}{x_k-x_j}\right)\frac{\partial u_m}{\partial x_j} .
\end{multline}
Putting together \eqref{der_v_xi}, \eqref{middle}, and \eqref{OmegaPum_xn}, \eqref{res_main}, we obtain \eqref{um_xkxk}.
$\Box$


\begin{corollary}
\label{cor_genus1}
Let $\T\to  X$ be a Toda family of genus one Riemann surfaces $\Tx$  corresponding to the elliptic curves of equation
\begin{equation}
\label{curve_g1}
\mu^2=\l(\l-1)(\l- x)(\l- u),
\end{equation}
parameterized by $\x=x$ as in Definition \ref{def_Todafamily}. Then the position of the  branch point $u$ of the coverings $\l:\Tx\to\mathbb CP^1$ as a function of $x$ satisfies the following equation
\begin{multline}
\label{iso_genus1}
u''=\frac{1}{2}\left( \frac{u}{x-1} - \frac{u-1}{x} + \frac{1}{u-x}  \right) -\frac{u'}{2} \left(\frac{2}{x}+\frac{2}{x-1}  + \frac{1}{u-x} \right)\\
+\frac{(u')^2}{2} \left( \frac{2}{u}+\frac{2}{u-1}+\frac{1}{x-u}\right) - \frac{(u')^3}{2} \left( \frac{x}{u-1}-\frac{x-1}{u}+\frac{1}{x-u}  \right)\,.
\end{multline}
\end{corollary}
{\it Proof.} This is the $g=1$ case of equation \eqref{um_xkxk} with $x_k=x$ and $u_m=u.$ Equations \eqref{um_xkxn} are empty in genus one. $\Box$

\begin{remark}
Note that, for a given initial value $u(x_0)$, the value of the derivative $u'=-\frac{ \Omega(P_x)\phi(P_x)}{\Omega(P_u)\phi(P_u)}$ with $\phi$ defined by \eqref{phi} singles out the unique solution of \eqref{iso_genus1} defining the Toda family of genus one curves \eqref{curve_g1}. This expression for $u'$ is the genus one case of \eqref{umder_g}. The existence of Toda family for real $x_0$ and $u(x_0)$ is given by Theorems 2.7 and 2.12 of \cite{PS1999}.
\end{remark}

\section{Generalized Toda families of Riemann surfaces}
\label{sect_generalized}


In this section we define a family of  hyperelliptic curves that generalizes the Toda family from Definition \ref{def_Todafamily}. All the results obtained  in the paper for a Toda family hold for these generalized families. Let us start by generalizing $\Omega$, the meromorphic differential of the third kind  \eqref{Omega} and define the following differential on the Riemann surface corresponding to the hyperelliptic curve \eqref{surf}
\begin{equation}
\label{Omegalpha}
\Omega_\alpha(P) = \Omega_{\infty^-\infty^+}(P) -\alpha^T\omega,
\end{equation}
where $\alpha\in\mathbb C^g$ is an arbitrary constant column vector, $\Omega_{\infty^-\infty^+}$ is the differential \eqref{OW} and $\omega=(\omega_1, \dots, \omega_g)^T$ is, as before, a vector of normalized holomorphic differentials from \eqref{B}. Setting $\alpha=2\pi\i N$ with $N$ from \eqref{Abel}, we recover definition \eqref{Omega} of $\Omega$ used in previous sections. The $a$-periods of $\Omega_\alpha$ are constant and given by components of the vector $\alpha.$ The $b$-periods are
\begin{equation}
\label{b-Omegalpha}
\oint_{\b_k}\Omega_\alpha =2\pi\i\int_{\infty^-}^{\infty^+}\omega_k - \alpha^T\mathbb B_k,
\end{equation}
where $\mathbb B_k$ is the $k$th column of the Riemann matrix $\mathbb B$ \eqref{B}.

Now, we would like to deform the hyperelliptic surface by varying the branch points of the corresponding two-fold covering in such a way that the $b$-periods \eqref{b-Omegalpha} stay constant. Assume we have such a family of surfaces $\Tx$ corresponding to the hyperelliptic curves of equation \eqref{surf} parameterized by the values of the branch points $x_1, \dots, x_g$ of the covering $u:\Tx\to\mathbb CP^1$ and having branch points $u_1, \dots, u_g$ as dependent functions. We assume that a meromorphic differential $\Omega_\alpha$ is defined on each $\Tx$ by \eqref{Omegalpha} and that its periods are constant for all $(x_1, \dots, x_g)$ varying continuously in some set. Let us then differentiate the right hand side of \eqref{b-Omegalpha} with respect to $x_i$, and equate the derivative to zero. With the help of the Rauch formulas \eqref{RauchB} for $\omega$ and $\mathbb B$, we have
\begin{multline*}
0 =\pi\i\,\omega_k(\pxi)\Omega_{\infty^-\infty^+}(\pxi)- \pi \i \,\omega_k(\pxi)\alpha^T\omega(\pxi)
\\
+\pi\i\,\sum_{j=1}^g\omega_k(\puj)\Omega_{\infty^-\infty^+}(\puj)\frac{\partial u_j}{\partial x_i}- \pi \i \sum_{j=1}^g\omega_k(\puj)\alpha^T\omega(\puj)\frac{\partial u_j}{\partial x_i}\,.
\end{multline*}
This is equivalent to the $k$th component of vector equation \eqref{forlater} with $\alpha=2\pi\i N$ and thus can be rewritten in the form \eqref{tempdiff_g} with $\Omega$ replaced by $\Omega_\alpha.$ This implies that Theorem \ref{thm_umder} giving an expression for $\partial u_m/\partial x_i$ in terms of $\Omega$ carries over to our new settings, for the isoperiodic deformations of $\Omega_\alpha$, upon replacing $\Omega$ by $\Omega_\alpha.$ Moreover, Theorem \ref{thm_support} also holds for the isoperiodic deformation as it is derived directly from Theorem \ref{thm_umder}, the Rauch formulas and the structure of the curves and differentials on them.

In other words, we may define a more general family of curves  including the Toda family of Definition \ref{def_Todafamily} as a particular case.
\begin{definition}
\label{def_Todatype}
We call  generalized Toda family or a  family of Toda type, a triple $(T, s_{\infty^+}, s_{\infty^-})$, if $ T: \T\to  X$ is a smooth fibration with fibers given by compact hyperelliptic genus $g$ Riemann surfaces  $\Tx$ for $\x\in X$ such that
\begin{itemize}
\item $X\subset\mathbb C^g\setminus \Delta$ with $\Delta=\{(x_1, \dots, x_g)\in\mathbb C^g\;|\; \exists\, k\neq j \mbox{ with } x_k=x_j \}$,
\item $u:\Tx\to\mathbb CP^1$ is a function of degree two on $\Tx$,
\item $\{\a_1,\dots, \a_g; \b_1, \dots, \b_g\}$ is a canonical homology basis on $\Tx$ for $\x\in X$ where $u(\a_k)$ and $u(\b_k)$ are independent of $\x\in X$ for all $1\leq k\leq g,$
\item $s_{\infty^+}$ and $s_{\infty^-}$ are two sections of $T$ such that  on every surface $\Tx$ of the family there exists a meromorphic differential of the third kind with simple poles at $s_{\infty^+}(\x)$ and $s_{\infty^-}(\x)$ with residues $+1$ and $-1$, respectively, whose $a$- and $b$-periods  are independent of $\x\in X.$
\end{itemize}
%

We also say that the generalized Toda family provides an isoperiodic deformation of the pair $(\T_\x, \Omega(\x))$, or an isoperiodic deformation of $\T_\x$ relative to $\Omega(\x)$, where $ \Omega(\x)$ is a meromorphic differential of the third kind on the surface $\T_\x$ for some fixed $\x$ having simple poles at $s_{\infty^+}(\x)$ and $s_{\infty^-}(\x)$ with residues $+1$ and $-1$.
\end{definition}
In this section, we showed that the following theorem holds.
\begin{theorem}
\label{thm_isoperiodicity}
Let $\T\to  X$ be a generalized Toda family of hyperelliptic Riemann surfaces $\Tx$  corresponding to the curves of equation \eqref{surf} parameterized by $\x=(x_1, \dots, x_g)\in X$ as in Definition \ref{def_Todatype}. Let $\Omega_\alpha$ \eqref{Omegalpha} be the differential of the third kind with poles at $s_{\infty^+}(\x)$ and $s_{\infty^-}(\x)$ with residues $+1$ and $-1$, respectively, from Definition \ref{def_Todatype}.  Then the  branch points $(u_1, \dots, u_g)$ of the coverings $u:\Tx\to\mathbb CP^1$ as functions of $\x=(x_1, \dots, x_g)$ have the derivatives expressed by
\begin{equation}
\label{umder_alpha}
\frac{\partial u_m}{\partial x_i}=-\frac{ \Omega_\alpha(\pxi)}{\Omega_\alpha(\pum)} v_m(\pxi)\,.
\end{equation}
and satisfy system \eqref{um_xkxn}, \eqref{um_xkxk} of Theorem \ref{thm_support}.
\end{theorem}
\begin{theorem}
\label{thm_converse-Toda}
Let $\x_0\in\mathbb C^g\setminus\{(x_1, \dots, x_g)\; |\; \exists\,i\neq j \text{  with }   x_{i}=x_j \}$ such that $x_{0j}\neq 0$ for $j=1, \dots, g$. Let $\Omega_\alpha=\Omega_\alpha(\x_0)$ be the meromorphic differential of the third kind defined by \eqref{Omegalpha} on a compact Riemann surface $\T_{\x_0}$ of the hyperelliptic curve \eqref{surf} of genus $g$ with $\x=\x_0$ and distinct $u_1, \dots, u_g,$ not coinciding with any of $x_1, \dots, x_g$.  Assume that $\x_0$ is such that $\Omega_\alpha(\x_0; \puj)\neq 0$ for $j=1, \dots, g$. Let $\{a_1, \dots, a_g,b_1, \dots, b_g\}$ be a canonical homology basis on $\T_{\x_0}$ such that the projections $u(a_j)$ and $u(b_j)$ with $j=1, \dots, g$ do not intersect a certain neighbourhood $\hat{\mathcal X}$ of $\x_0$. Then there exists a unique continuous  $g$-parameter generalized Toda family  providing an isoperiodic deformation of the pair $(\T_{\x_0}, \Omega_\alpha(\x_0))$  parametrized by $\x$ varying in some neighbourhood $\mathcal X\subset\hat{\mathcal X}$ of $\x_0.$
\end{theorem}
{\it Proof.}
The isoperiodic deformations of $(\T_{\x_0}, \Omega_\alpha(\x_0))$ are defined by $g$ conditions  $\beta_k=\oint_{b_k}\Omega_\alpha(x)={\rm const},$ $k=1, \dots, g$.  These relations define implicitly the functions $\u:=(u_1(\x), \dots, u_g(\x))$ in a neighbourhood of $\x=\x_0$ if the the Jacobian ${\rm det} J_\u(\x,\u) $ is non-zero at $\x=\x_0$, where the $g\times g$ matrix $J_\u(\x,\u)$ has the following $jk$-entry $J^{jk}_\u(\x,\u)$:
\begin{equation*}
 J^{jk}_\u(\x,\u)= \frac{\partial}{\partial u_j} \oint_{b_k}\Omega_\alpha(\x_0).
\end{equation*}
From \eqref{b-Omegalpha} and \eqref{Omegalpha} we obtain these derivatives using the Rauch formulas \eqref{RauchB} and the relation \eqref{OW}:
\begin{equation*}
J^{jk}_\u(\x,\u) = 2\pi\i\int_{\infty^-}^{\infty^+}\frac{\partial^{\rm Rauch}}{\partial u_j}\omega_k - \alpha^T \frac{\partial^{\rm Rauch}\mathbb B_k}{\partial u_j}=
 \pi\i\,\Omega_\alpha(\x;\puj)\omega_k(\x;\puj).
\end{equation*}
Thus ${\rm det} J_\u(\x_0,\u) = \pi\i\,{\rm det} J \prod_{m=1}^g\Omega_\alpha(\x_0;\pum)$ where by $ J$ we denote the matrix $ J_{jk} = \omega_k(\x_0;\puj).$ The product of $\Omega_\alpha(\pum)$ being non-zero at $\x=\x_0$, we only need to show that ${\rm det} J\neq 0.$ Note that the vanishing of ${\rm det} J$ is equivalent to the vanishing of ${\rm det} \hat J$ where $\hat J$ is a matrix constructed similarly to $J$ from another basis $\phi(P), \;\lambda(P)\phi(P), \dots, \lambda(P)^{g-1}\phi(P)$ \eqref{phi} in the space of holomorphic differentials. More precisely, $\hat J_{jk} = u_j^{k-1}\phi(\x_0;\puj)$.
The determinant of $\hat J$ is the Vandermonde determinant $V(u_1, \dots, u_g)$, which is non-zero for distinct $u_1, \dots, u_g$. This being fulfilled at $\x=\x_0$, we proved the existence of local continuous isoperiodic deformation of $(\T_{\x_0}, \Omega_\alpha(\x_0))$ and thus of the generalized Toda family.

Now, Theorem \ref{thm_isoperiodicity} shows that for every continuous isoperiodic deformation, the functions $u_1(\x), \dots, u_g(\x)$ satisfy equations \eqref{um_xkxn} and \eqref{um_xkxk} for all $1\leqslant m,n,k\leqslant g$. In addition,  the first derivatives $\partial_{x_k}u_j(\x)$ are given by \eqref{umder_alpha} for any $1\leqslant k,j\leqslant g$ and $\x=(x_1,\dots, x_g)$. Thus, for $\x_0=(x_1^0, \dots, x_g^0)$, the  initial values $u_j(\x_0)$ and first order derivatives \eqref{umder_alpha} for any $\x$ specify a unique solution of system \eqref{um_xkxn}, \eqref{um_xkxk}  as follows. On the set of points $(x_1, x_2^0 \dots, x_g^0)$ for $x_1$ varying in some neighbourhood of $x_1^0$, the system reduces to a system of ordinary differential equations of second order with respect to the variable $x_1$, and thus the initial values $u_j(\x_0)$ and the derivatives $\partial_{x_1}u_j(\x_0)$ single out a unique solution $\u(x_1, x_2^0 \dots, x_g^0)$ for  $x_1\in S_1$ for some open neighbourhood $S_1$ of $x_1^0$. Now, the obtained functions $\u(x_1, x_2^0 \dots, x_g^0)$ give us initial values for the system of ordinary differential equations of second order with respect to the variable $x_2$ at any point of the set $\{(x_1, x_2^0 \dots, x_g^0)|\; x_1\in S_1\}$. Together with the values of derivatives $\partial_{x_2}u_j((x_1, x_2^0 \dots, x_g^0))$, they single out a unique solution $\u(x_1, x_2, x_3^0 \dots, x_g^0)$ for  $(x_1,x_2)\in S_2$ for some open neighbourhood $S_2$ of $(x_1^0, x_2^0)$. Proceeding in this way, we obtain the uniqueness of an isoperiodic deformation of $(\T_{\x_0}, \Omega_\alpha(\x_0))$ locally in an open neighbourhood of $\x_0$.
$\Box$

\section{Isoequilibrium deformations}
\label{sect_isoeq}

For a measure $\mu$ defined on $\mathbb C$, one defines its {\it potential energy} by
\begin{equation*}
\mathcal E(\mu)=\int_{\mathbb C^2} \log (|z-w|^{-1})|d\mu(z)d\mu(w).
\end{equation*}
Let $c\subseteq \mathbb C$ and consider all probability measures for which the support is a subset of $c$. If the energy of all such probability measures
is equal to infinity, the set $c$ is called of {\it capacity zero.} If  $ c$ is not of capacity zero, then there exists  a unique probability measure  $\rho_{e}(c)$ whose support is a subset of $c$ and whose energy is minimal among all probability measures supported in subsets of $c$,  see e.g. \cite{Si2011}, \cite{Si2015a}. The measure  $\rho_{e}(c)$ is called {\it the equilibrium measure} of $c$.

Consider $c$ to be a union of $g+1$ nonintersecting   intervals
%
$c=\bigcup_{j=0}^{g}[\alpha_j, \beta_j].$
%
The Lebesgue measure restricted to $c$ and appropriately scaled has a finite potential energy, therefore $c$ is not of capacity zero. By using  an affine transformation, one can map $c$ into a union of the interval $[0,1]$ and  $g$ intervals with the end-points denoted by $x_j, u_j$, $j=1, \dots, g$.
Let us denote the M\"obius transformed set of intervals also by $c$, that is we have
\begin{equation}
\label{c}
c=[0,1]\cup\bigcup_{j=1}^{g}[x_j, u_j], \quad \text{assuming}\quad 1<x_j<u_j<x_{j+1}.
\end{equation}
Associated with this set of real intervals, there is a compact hyperelliptic Riemann surface $\surf_\x$ corresponding to the algebraic curve defined by \eqref{surf}, or by \eqref{T} and \eqref{delta}. As before, the points $0,1,x_j, u_j$ are the branch points of the covering $u:\surf_\x\to\mathbb CP^1$ with $\x=(x_1, \dots, x_g)$  and $P_0, P_1, \pxj, \,\puj$ are the corresponding ramification points on the surface $\surf_\x$.

  Let $\{\a_1,\dots, \a_g; \b_1, \dots, \b_g\}$ be a fixed canonical homology basis on $\surf_\x,$ such that $u(\a_1)$ goes around the branch points 1 and $x_1$, and $u(\a_j)$ goes around the  points $u_{j-1}$ and  $x_j$; while $u(\b_1)$ goes  clockwise around the first interval $[0,1]$, and $u(\b_j)$ goes  clockwise around the first $j$ intervals, from $[0,1]$ to $[x_{j-1}, u_{j-1}]$. The $a$-cycles are oriented in a way to obtain the following intersection indices: $\a_k\circ \a_j=\b_k\circ \b_j=0$ and $\a_k\circ \b_j=\delta_{kj}$.

As in Section \ref{sect_TodaPell}, there are two points at infinity on $\surf_\x$, denoted by $\infty^+$ and $\infty^-$, and we may define a differential $\Omega_{\alpha}$  by \eqref{Omegalpha}.
We consider a special case  of $\alpha=0$, that is $\Omega_0=\Omega_0(\x)=\Omega_{\infty^-\infty^+}$ is a differential \eqref{OW} of the third kind with poles at $\infty^\pm$ with residues $\pm 1$ normalized by vanishing $a$-periods. According to \cite{Widom}, \cite{Si2011}, we have for the equilibrium measures  (see e.g. equation (5.6.4) in \cite{Si2011} and Section 2.2 in \cite{DR2023}):
\begin{equation}
\label{eqmeasure}
 \rho_{e} ([x_j, u_j])=  \frac{1}{\pi}\int_{x_j}^{u_j}\frac{|q_{g}(x)|d\,x}{\sqrt{|\Delta_{2g+2}(x)|}}\quad j=1,\dots, g, \quad \text{assuming}\quad x_j<u_j.
\end{equation}
 Here  $\Delta_{2g+2}$ is given by \eqref{delta} and $q_{g}$ is the monic polynomial of degree $g$ defined by the conditions:
\begin{equation}
\label{null}
\int_{u_{j-1}}^{x_j}\frac{q_{g}(x)d\,x}{\sqrt{\Delta_{2g+2}(x)}}=0, \quad j=1, \dots, g,
\end{equation}
where $u_{0}=1$.

\begin{lemma} The equilibrium measure is related to the $b$-periods of the differential $\Omega_0$ in the following way:
\begin{equation}
\label{eq:bjrho}
\oint_{\b_1}\Omega_0=2\pi\i\rho_{e} ([0, 1]);  \oint_{\b_2}\Omega_0= 2\pi\i(\rho_{e} ([0, 1])+\rho_{e} ([x_1, u_1])); \ldots, \oint_{\b_g}\Omega_0=2\pi\i\big(\rho_{e} ([0, 1])+\sum_{j=1}^{g-1}\rho_{e} ([x_{j}, u_{j}])\big),
\end{equation}
where
$$
 \rho_{e} ([0, 1])=1-\sum_{j=1}^g\rho_{e} ([x_j, u_j]).
$$
\end{lemma}

{\it Proof.}  From \eqref{null}, we see that $q_{g}$ changes sign on each interval of the form $(u_{j-1}, x_j)$. Since the polynomial  $q_{g}$  is of degree $g$ and in each of $g$ intervals $(u_{j-1}, x_j)$ it has at least one zero, it has exactly one zero in each of these intervals. It has no zeros outside these intervals. Thus, it has a constant sign on each of the intervals $[x_j, u_j]$ and this sign alternates on consecutive intervals of this form.  On the Riemann surface $\T_\x,$ from the assumption that all $a$-periods of $\Omega_0$ are zero  and from the structure of its poles, we get
\begin{equation}
\label{Omega-equilibrium}
\Omega_0=-\frac{q_{g}(u) du}{\sqrt{\Delta_{2g+2}(u)}}.
\end{equation}

 Assume now that the branch cuts of the covering $u:\surf_\x\to\mathbb CP^1$ coincide with the intervals composing the set $c$ \eqref{c}.
 The sign of the analytic continuation of $\i\sqrt{\Delta_{2g+2}(x)}$ alternates on the intervals $[x_j, u_j]$ and the branch is chosen so that
 $\i\sqrt{\Delta_{2g+2}(x)}$ is negative on the upper side of the cuts  along those intervals where $q_{g}$ is positive. Thus, we get:

 \begin{equation*}
 \frac{1}{\pi}\int_{x_j}^{u_j}\frac{|q_g(x)|d\,x}{\sqrt{|\Delta_{2g+2}(x)|}}    =\frac{1}{2\pi\i}\oint_{\hat b_j}\Omega_0,\quad j=1,\dots, g,
 \end{equation*}
 where $\hat b_j$ is a cycle on $\T_\x$ such that $u(\hat b_j)$ goes clockwise once around $[x_j, u_j]$, that is $\hat b_j=\bf b_{j+1}-\bf b_{j}$ for $j=1, \dots, g.$ This is due to the fact that $\Omega_0$ satisfies $\Omega_0(\x;P)=-\Omega_0(\x;P^*)$, where $P^*\in\T_\x$ denotes the hyperelliptic involution applied to $P\in\T_\x.$
 Comparing to \eqref{eqmeasure},  we see that $ \rho_{e} ([x_j, u_j])$, for $j=1, \dots, g$ uniquely determine the $b$-periods of $\Omega_0$ and vice-versa, and we get
  \eqref{eq:bjrho}.
$\Box$
\\

 Thus, the deformations of $x_j$, $ u_j$  preserve  $ \rho_{e} ([x_j, u_j])$, for $j=1, \dots, g$ if and only if they preserve the $b$-periods of $\Omega_0$.
Our goal now is to study smooth variations of the endpoints of the intervals in a way that the equilibrium measures remain constant.
\begin{definition}
\label{def_isoequilibrium}
Let $c$ be a set of real intervals given by \eqref{c}. An isoequilibrium deformation of $c$ is a smooth transformation of the intervals $[x_j, u_j]$, $j=1, \dots, g$,  along the real line, which fixes the first interval $[0,1]$  of $c$,  and which leaves
 the equilibrium measures of the intervals $[x_j, u_j]$ for $j=1, \dots, g$ unchanged.
\end{definition}

\begin{proposition}
\label{prop:isoequirat}
Given a set of real  intervals $c(\x)$ \eqref{c} such that  the equilibrium measure of each interval  is rational,
 there exists a unique  $g$-parameter isoequilibrium deformation  of $c(\x)$  produced by varying $\x$ in some neighbourhood  $\mathcal X$ in $\mathbb R^g$ and having $u_1, \dots, u_g$ as real-valued functions of $\x$   for   $\x\in\mathcal X$.
\end{proposition}

{\it Proof.} It was shown in \cite{Apt1986} that the  for a set of real  intervals $c(\x)$ \eqref{c},  the equilibrium measure of each interval is rational if and only if the set of real  intervals \eqref{c} is a support of a generalized Chebyshev polynomial. The rational numbers describing the equilibrium measures translate into the combinatorial properties of the polynomial (see Remark \ref{rem:rat} below). From \cite{PS1999}, Theorems 2.7 and 2.12,  follows the existence and uniqueness of a  continuous deformation of the  generalized Chebyshev polynomial preserving its combinatorial structure, so that its support is a continuous deformation of the set of intervals \eqref{c}, which preserves the isoequilibrium measures.
$\Box$
\\

Isoequilibrium deformations of the set of intervals induce deformations of the corresponding surface $\mathcal T_\x$. In the next corollary we show that these deformations of surfaces are the isoperiodic deformations of $\mathcal T_\x$ in the sense of Definition \ref{def_Todatype}  relative to the differential $\Omega_0(\x)$ defined above.
\begin{corollary}
\label{cor:isoequi}
Let $\x=(x_1, \dots, x_g)$ and $c=c(\x)$ be a set of  real intervals given by \eqref{c}. Let $\T_\x$  be the associated hyperelliptic Riemann surface corresponding to the algebraic curve of equation \eqref{surf}. Let $\pxj$ and $\puj$ be the ramification points of the covering $u:\Tx\to\mathbb CP^1$ such that $u(\pxj)=x_j$ and $u(\puj)=u_j$ for $j=1, \dots, g.$ Let $c(\x)$ be subject to
isoequilibrium deformation parametrized by varying $\x$ in some neighbourhood  $\mathcal X$  in $\mathbb R^g$. Then  and only then the associated surfaces $\T_\x$  with the canonical homology basis chosen as above in this section and differentials $\Omega_0(\x)$ defined by \eqref{Omegalpha} with $\alpha=0$ on $\T_\x$  form a generalized Toda family,  such that the  branch points $u_1, \dots, u_g$ of the coverings $u:\Tx\to\mathbb CP^1$ as functions of $\x$ are real and have derivatives expressed by \eqref{umder_alpha}, \eqref{v} with $\alpha=0$
and satisfy system \eqref{um_xkxn}, \eqref{um_xkxk} of Theorem \ref{thm_support}.
\end{corollary}
{\it Proof.} Let $c(\x)$ for $\x$ belonging to some neighbourhood $\mathcal X\subset \mathbb R^g$ be an isoequilibrium  deformation of the set of intervals \eqref{c}.
 Let us prove that the corresponding pairs $(\T_\x, \Omega_0(\x))$ form a generalized Toda family.  The $a$-periods of $\Omega_0(\x)$ are zero by definition \eqref{Omegalpha} since $\alpha=0$. The $b$-periods  uniquely determine  the integrals in \eqref{eqmeasure}, and thus they uniquely determine the equlibrium measures for the real intervals and vice-versa, see \eqref{eq:bjrho}.
Definition \ref{def_isoequilibrium} implies that the integrals \eqref{eqmeasure} are independent of $\x\in\mathcal X$, and therefore the $b$-periods of $\Omega(\x)$ are independent of $\x$ as well and vice-versa. Thus, pairs $(\T_\x, \Omega_0(\x))$ form a generalized Toda family in the sense of Definition \ref{def_Todatype}  with the functions  $u_j(\x)$ being real.
The  statement concerning differential equations satisfied by $u_j(\x)$ follows from Theorems \ref{thm_support} and \ref{thm_isoperiodicity}.
$\Box$

\begin{proposition}\label{prop:isoequiirat}
 Given a set of real  intervals $c(\x)$ \eqref{c} with arbitrary equilibrium measures,
 there exists a unique  $g$-parameter isoequilibrium deformation  produced by varying $\x$ in some neighbourhood  $\mathcal X$ in $\mathbb R^g$ and having $u_1, \dots, u_g$ as real-valued functions of $\x$.
\end{proposition}

{\it Proof.}  Consider the hyperelliptic curve $\mathcal T_\x$ \eqref{surf} associated with the set of real intervals $c(\x)$ with the homology basis and differential $\Omega_0$ as in Corollary \ref{cor:isoequi}. Assume the pair $(\mathcal T_\x, \Omega_0(\x))$ is a generalized Toda family for some $\x\in\hat{\mathcal  X}.$ We want to prove that for $\x$ varying in some subset of $\hat{\mathcal  X}\cap \mathbb R^g$, the functions $u_1(\x), \dots, u_g(\x)$ are real.

According to  the Bogatyrev--Peherstoefer--Totik Theorem, see Theorem 5.6.1 in \cite{Si2011} and its proof, we know  that the set of real  intervals \eqref{c} can be approximated by sets of intervals,
\begin{equation}
\label{ckn}
c^{(n)}=[0,1]\cup\bigcup_{j=1}^{g}[x^{(n)}_j, u^{(n)}_j], \quad 1<x^{(n)}_j<u^{(n)}_j<x^{(n)}_{j+1},
\end{equation}

such that each interval has a rational equilibrium measure and such that the difference between $c^{(n)}$ and $c$ tends to zero as $n$ tends to infinity, that is $x^{(n)}_j\to x_j$ and $u^{(n)}_j\to u_j$. According to Proposition \ref{prop:isoequirat}, each of the sets $c^{(n)}$ admits an isoequilibrium deformation parametrized by ${\bf x}^{(n)}=(x_1^{(n)}, \dots, x_g^{(n)})\in\mathcal X\subset\mathbb R^g$. According to Corollary \ref{cor:isoequi}, such isoequlibrium deformations generate  isoperiodic deformations of the associated hyperelliptic Riemann surface $\mathcal T_{{\bf x}^{(n)}}$ given by \eqref{surf} with $\bf x$ replaced by ${\bf x}^{(n)}$ relative to the differential $\Omega_0$ given by \eqref{Omega} with $\alpha=0$ on $\mathcal T_{{\bf x}^{(n)}}$ and appropriately chosen homology basis.  Moreover, the  functions $u^{(n)}_1({\bf x}^{(n)}), \dots, u^{(n)}_g({\bf x}^{(n)})$ are real under these deformations. Thus, their derivatives expressed according to \eqref{umder_alpha} are real for each $c^{(n)}$.

In the limit, when $n$ tends to infinity, the surfaces $\mathcal T_{{\bf x}^{(n)}}$ tend to $\mathcal T_\x$ and the chosen of homology bases on $\mathcal T_{{\bf x}^{(n)}}$ produce the required homology basis on $\mathcal T_\x$. Thus the differentials $\Omega$ and $v_j$ on $\mathcal T_{{\bf x}^{(n)}}$ become those on $\mathcal T_\x$. Therefore, the derivatives \eqref{umder_alpha} for the functions $u^{(n)}_1({\bf x}^{(n)}), \dots, u^{(n)}_g({\bf x}^{(n)})$ tend to derivatives of the functions $u_1(\x), \dots, u_g(\x)$ with $\x$ restricted to some real subset of $\hat{ \mathcal X} \subset \mathbb R^g$. We thus may conclude that  derivatives $\partial u_j/\partial x_k$ are real. The functions $u_j$ being solutions of the real system \eqref{um_xkxn}, \eqref{um_xkxk} of Theorem \ref{thm_support} with real initial values are then real as well.
Thus the restriction of the original Toda family $(\mathcal T_\x, \Omega_0(\x))$ on some real subset of $\hat{\mathcal X}$ for $\x$ has $u_1,\dots, u_g$ real and thus corresponds to  an isoequilibrium deformation of $c(\x)$.
$\Box$

\subsection{Deformations of the essential spectra of periodic Jacobi matrices}

Following \cite{Si2011}, we now consider infinite two-sided  tridiagonal Jacobi matrices $J$ determined by two two-sided sequences $\{(c_n), (d_n)\}_{n=-\infty}^{+\infty}$. The Jacobi matrices are defined as operators acting on $\ell^2(\mathbb Z)$  in the following way: for $f=(f_n)_{n=-\infty}^\infty\in\ell^2(\mathbb Z)$ the elements of the sequence $Jf$ are given by
\begin{equation*}
(Jf)_n=c_nf_{n+1}+d_nf_n+c_{n-1}f_{n-1}.
\end{equation*}
The Jacobi matrix $J$ is periodic with period $N$ if the sequences $c_n, d_n$ are both periodic with period $N$.
Generically, the essential spectrum of an $N$-periodic Jacobi matrix consists of $N$  real  intervals; the $N-1$ intervals between them are called {\it gaps}. It may happen that some of the intervals merge eliminating some of the gaps, which results in the essential spectrum consisting of $g+1\leqslant N$  intervals.
 From \cite{Per}  and \cite{Apt1986} (see also e.g. \cite{Si2011}, Theorem 5.13.8) we have the following statement.

\begin{theorem}
\label{cor:rational}
A union of intervals $c$ \eqref{c} is a support of  a generalized Chebyshev polynomial of degree $N$ if and only if it is
the essential spectrum of a periodic Jacobi matrix $J$ with period $N$. This is the case if and only if all the measures   $\rho_{e}([x_j, u_j])$  are rational and  $N\rho_{e}([x_j, u_j])$   are integer for all $j$.
\end{theorem}

{\it Proof.}
Let a set $c$ of the form \eqref{c} be an essential spectrum of an $N$-periodic Jacobi matrix $J$.
For an essential spectrum $c$ of $J$, by formula (5.1.14) in \cite{Si2011},  there exist  natural numbers   $k_j$,  $j=0, \dots, g$,  such that
\begin{equation*}
 \rho_{e}([0, 1])=\frac{k_0}{N}, \quad \rho_{e}([x_j, u_j])=\frac{k_j}{N}, \quad  j=1, \dots, g, \quad \sum_{j=0}^gk_j=N.
 \end{equation*}

For a union of real intervals, there exists a natural number $N$ such that the equilibrium measure of each interval is a rational number, which multiplied by $N$ gives an integer, if and only if there
exists a conformal map from the upper half plane to a comb region $\mathcal C_N$  of the form
$$
\mathcal C_N=\{w|\, \Im w>0, 0\le \Re w\le N\pi\}\setminus\bigcup_{j=1}^{N}\{w|\, \Re w = j\pi, 0\le \Im w\le h_j\},
$$
with some nonnegative values  $h_1, \dots, h_{N}$. In other words,
this comb region is a vertical semi-strip with a finite number of  vertical slits where the {\it basis} of the comb region is a real segment together with the marked points  $j\pi$ on it. The conformal map maps the union of the intervals onto the basis of the comb and $\infty$ to $\infty$, see e.g. \cite{SY1992}, Section 5.

According to \cite{Per}, Theorem 1, such a conformal map exists for a union of real intervals with $N\in\mathbb N$ if and only if it is a support of the Pell equation with a Chebyshev polynomial of degree $N$, see also e.g. \cite{SY1992}, Theorem of Section 5.5.

For detailed proofs that every support of a Pell's equation with a polynomial of degree $N$ is the essential spectrum of a periodic  Jacobi operator with the period $N$, see  \cite{Per} and also e.g. \cite{Si2011}, Theorem 5.13.8.
$\Box$\\

 For the hyperelliptic surface associated with a set of intervals to be a Toda curve, we need the surface  to satisfy \eqref{Abel} with $M_1, M_2$ such that $NM_1, NM_2\in\mathbb Z^g$. This is a weaker condition than \eqref{be-qmeasure} due to the presence of a possibly nonzero $M_2$ in \eqref{Abel}.
\begin{corollary}\label{cor:M20}  Given a Toda curve of order $N$ in the sense of Definition \ref{def_Toda}, which satisfies \eqref{Abel} with $M_1, M_2$ such that $NM_1, NM_2\in\mathbb Z^g$, the condition
$M_2=0$ in \eqref{Abel} is equivalent to the condition that there exists a generalized Chebyshev polynomial of degree $N$ that satisfies the Pell equation on the set of intervals $[0,1]$, $[x_j, u_j]$, for $j=1,\dots, g$.
\end{corollary}
{\it Proof.}  Let us start from a surface $\T_\x$, which is a Toda curve of order $N$ according to Definition \ref{def_Toda}.
We denote  $\x=(x_1, \dots, x_g)$ and assume the choice of a canonical homology basis as above  in Section \ref{sect_isoeq}. We define $\Omega_0(\x)$ by \eqref{Omegalpha} with $\alpha=0$ as above as well. The equilibrium measures of the intervals $[0,1]$, $[x_j, u_j]$, for $j=1,\dots, g$ are of the form

\begin{equation*}
 \rho_{e}([0, 1])=\frac{k_0}{N}, \quad \rho_{e}([x_j, u_j])=\frac{k_j}{N}, \quad k_j\in\mathbb N, \quad  j=1, \dots, g, \quad \sum_{j=0}^gk_j=N,
 \end{equation*}
if and only if,  by \eqref{eqmeasure} and \eqref{eq:bjrho},   the $b$-periods of $\Omega_0(\x)$ are rational numbers up to a factor of $2\pi\i$ of the following form:
\begin{equation}
\label{be-qmeasure}
\oint_{b_j}\Omega_0 =\frac{2\pi\i \sum_{\ell=0}^{j-1}k_{\ell}}{N} ,\quad j=1,\dots, g.
\end{equation}
Thus, according to \eqref{b-Omegalpha}, $M_2=0$ in \eqref{Abel} if and only if the equilibrium measures of the intervals $[0,1]$, $[x_j, u_j]$, for $j=1,\dots, g$ are all rational.

The proof now follows from Theorem \ref{cor:rational} and its proof and \cite{Per}, Theorem 1.
$\Box$

\begin{remark}\label{rem:rat}
The numbers $k_j$,  $j=0, \dots, g$  from \eqref{be-qmeasure} are determined by the structure  of the merging of the intervals, leading to the elimination of some of the  gaps.  From e.g. \cite{Si2011} one can conclude that the number  $k_0$  gives the number of intervals that  merged to form the leftmost interval $[0,1]$, $k_1$ is the number of the intervals that  merged to form the second interval from left $[x_1,u_1]$, and so on.
\end{remark}

\begin{corollary}\label{cor:Jacobiperiodic}  Let $\x_0=(x_1^0, \dots, x^0_g)$ and $c_0=c(\x_0)$ be a union of $g+1$  real intervals given by \eqref{c} with $\x=\x_0$. Let $c_0$ be the essential spectrum of an $N$-periodic Jacobi matrix. The sets of intervals $c=c(\x)$ \eqref{c} for $\x=(x_1, \dots, x_g)\in\mathcal X\subset\mathbb R^g$ obtained from $c_0$ by local smooth deformations of $\x$ that keep $[0,1]$  unchanged are the essential spectra of $N$-periodic Jacobi matrices with the same structure of gap intervals if and only if they are isoequlibrium deformations   parametrized by $\x$.  The associated hyperelliptic Riemann surfaces $\T_\x$ \eqref{T}, \eqref{delta} form a Toda family    of order $N$, with  $\x\in \mathcal X\subset \mathbb R^g,$  relative to the differential $\Omega_0$ given by \eqref{Omega-equilibrium}. The  branch points $(u_1, \dots, u_g)$ of the coverings $u:\Tx\to\mathbb CP^1$ are real-valued functions of $\x=(x_1, \dots, x_g)\in \mathcal X\subset \mathbb R^g$, have derivatives expressed by
\eqref{umder_alpha}, \eqref{v} with $\alpha=0$, and satisfy system \eqref{um_xkxn}, \eqref{um_xkxk} of Theorem \ref{thm_support}.
\end{corollary}

{\it Proof.} According to Theorem \ref{cor:rational}, a set of intervals $c_0=c(\x_0)$ is the essential spectrum of an $N$-periodic Jacobi matrix if and only if all the equlibrium measures are rational. According to Proposition \ref{prop:isoequirat},  there exists a unique  $g$-parameter isoequilibrium deformation over a small neighborhood $\x=(x_1, \dots, x_g)\in\mathcal X\subset\mathbb R^g$. Due to Corollary \ref{cor:isoequi}, this isoequilibrium deformation coincides with the restriction to $\x\in \mathcal X$ of the isoperiodic deformation of the associated hyperelliptic surface $\mathcal T_\x$ relative to the differential $\Omega_0$ from Corollary \ref{cor:isoequi}.  Thus, the  branch points $(u_1, \dots, u_g)$ of the coverings $u:\Tx\to\mathbb CP^1$ are real-valued functions of $\x=(x_1, \dots, x_g)\in \mathcal X\subset \mathbb R^g$ and they have derivatives expressed by
\eqref{umder_alpha}, \eqref{v} with $\alpha=0$, and satisfy system \eqref{um_xkxn}, \eqref{um_xkxk} of Theorem \ref{thm_support}.
$\Box$

  \begin{remark}\label{ex:Peiers} Isoequilibrium deformations both for the general case of arbitrary periods \eqref{eqmeasure} and in the special case where all periods are rational were studied in the physics oriented setting of Peierls models in \cite{BDK}. Such transformations were called  there ``variations that preserve the group of periods of the differential $dp$". In \cite{BDK}, extremum equations with respect to variations that preserve the group of periods were derived for the functionals of the total energy under the consideration in \cite{BDK}.
  \end{remark}

\section{Isoperiodic deformations of periodic Toda lattices}
\label{sect_lattices}
Toda lattice was introduced more than half a century ago \cite{To} and since then it occupies one of the central places in the theory
of integrable systems, see e.g \cite{Fla1, Fla2, Moer, FlaMcL, Nov, HK, Si2011} and references therein.

Algebro-geometric approach to Toda lattice has a long and elaborated history. One of the first steps were made by S. P. Novikov
in an unpublished paper which was incorporated into \cite{DMN}. These ideas were further developed by I. M. Krichever and others, see e.g. \cite{Krich1978, Nov, Krich1981} and references therein. Here, we will mostly follow \cite{Krich1981}, which appeared as an appendix of \cite{Dub1981}.

The Toda lattice describes a one-dimensional chain of particles with exponential interaction of immediate neighbours
\begin{equation}\label{eq:toda0}
\ddot x_n(t)=e^{x_{n+1}-x_n}-e^{x_{n}-x_{n-1}},
\end{equation}
where $x_n(t)$ is the position of the $n$-th particle at the moment $t$, and $\dot x_n(t)$ denotes the derivative with respect to $t$. Denoting $v_n:=\dot x_n$ and $c_n=\exp(x_n-x_{n-1})$, $c_n>0$,  Toda lattice equations \eqref{eq:toda0} can be rewritten in the form
\begin{equation}\label{eq:toda}
\begin{aligned}
\dot v_n&=c_{n+1}-c_n,\\
\dot c_n&=c_n(v_n-v_{n-1}).
\end{aligned}
\end{equation}
The integration procedure from \cite{Krich1981} starts from a hyperelliptic curve of the equation
\begin{equation}\label{eq:Gamma}
 \mu^2=\prod_{i=1}^{2g+2}(z-z_i)
\end{equation}
and the corresponding compact Riemann surface $\Gamma$ of genus $g$ with two points at infinity, $\infty^+$ and $\infty^-$. The central role in the integration is played by Krichever's Baker-Akhiezer
function $\psi(n, t, P)$ which is meromorphic on $\Gamma \setminus \{\infty^+, \infty^-\}$ having $g$ poles and zeros away from the two points at infinity and an essential singularity at $\infty^+, \infty^-$ with the following asymptotics
\begin{equation}
\label{eq:BA}
\psi(n, t, P)\underset{P\rightarrow\infty^{\pm}}{=}\i^n\lambda_n^{\pm 1}z^{\pm n}\left(1+\xi_1^{\pm}(n, t)z^{-1}+\dots\right)e^{\mp\frac{tz}{2}},
\end{equation}
where $\lambda_n, \,\xi_1^{\pm}(n, t)\in\mathbb C$.
For the difference operators $\mathbb L=(\mathbb L^{nm})$ and $\mathbb A=(\mathbb A^{nm})$, defined as follows:
\begin{equation}\label{eq:L}
\mathbb L^{nm}=-\i\sqrt{c_{n+1}}\delta_{n,m-1}+v_n\delta_{n,m}+\i\sqrt{c_{n}}\delta_{n,m+1},
\end{equation}
\begin{equation}\label{eq:A}
\mathbb A^{nm}=\frac{\i}{2}\sqrt{c_{n+1}}\delta_{n,m-1}+w_n\delta_{n,m}+\frac{\i}{2}\sqrt{c_{n}}\delta_{n,m+1},
\end{equation}
with $2(w_n-w_{n-1})=v_n-v_{n-1}-\ln c_n$,
we have
\begin{equation}\label{eq:LA}
\dot \psi=\mathbb A \psi, \quad \mathbb L\psi=z\psi.
\end{equation}
The Toda lattice equations \eqref{eq:toda} are the compatibility conditions for \eqref{eq:LA}. Using this, solutions to the Toda equations are obtained in terms of the coefficients of the expansion \eqref{eq:BA} for the Baker-Akhiezer function:  $c_n=(\lambda_{n-1}/\lambda_n)^2$, and $v_n=\xi^+(n+1,t)-\xi^+(n, t)$.
One can express the Baker-Akhiezer function $\psi(n, t, P)$ using the theta function of the surface $\Gamma$ corresponding to the curve \eqref{eq:Gamma} (see \cite{Dub1981} for the notation and properties of theta functions noting that the normalization in \cite{Dub1981} is different from the one we adopted in Section \ref{sect_differentials} )  and calculate $\lambda_n$ and $\xi^+_1(n, t)$ in \eqref{eq:BA} in terms of the theta functions as well. This implies the following expressions for  the solutions
of the Toda lattice \eqref{eq:toda}:
\begin{equation}
\label{eq:solTodav}
v_n=\frac{d}{dt}\ln\frac{\theta((n+1)U+tV+z_0)}{\theta(nU+tV+z_0)};
\end{equation}
\begin{equation}
\label{eq:solTodac}
c_n=\frac{\theta((n+1)U+tV+z_0)\theta((n-1)U+tV+z_0)}{\theta^2(nU+tV+z_0)},
\end{equation}
where $z_0$ is an arbitrary vector.
The vector $U$ is given by
\begin{equation}
\label{eq:vectorU}
U=2\pi\i \int_{\infty^-}^{\infty^+}\omega
\end{equation}
with $\omega=(\omega_1, \dots, \omega_g)^T$ being the vector of linearly independent holomorphic differentials on $\Gamma$ \eqref{eq:Gamma} normalized by the condition $\oint_{\a_j}\omega_k=\delta_{jk}$ (note that the normalization in \cite{Krich1981} is different).
The vector $V$ is a linear combination of the $b$-periods  of the normalized differentials of the second kind $\Omega_{\infty^-}$ and $\Omega_{\infty^+}$
having only a double pole at $\infty^-$ and $\infty^+$, respectively, and no other singularities.

As observed many times before (see e.g. \cite{Nov}), for the famous continuous integrable PDEs, like Korteweg-de Vries equation or the non-linear Schr\"odinger equation, it is practically impossible to effectively select solutions periodic in $x$ with a given period $T$. For the Toda lattice as a discrete system, the situation is remarkably different.
From \cite{Krich1981} we can extract the following
\begin{theorem} (\cite{Krich1981})
\label{th:todaN}
A solution to the Toda lattice \eqref{eq:toda} is periodic in $n$ with a period $N$ if and only if the solution $(v_n, c_n)$ is given by \eqref{eq:solTodav} and \eqref{eq:solTodac} relative to a hyperelliptic curve $\Gamma$ \eqref{eq:Gamma}  and the vector $U/(2\pi\i)$ given by \eqref{eq:vectorU} is an $N$-division point of the lattice generating the Jacobian of $\Gamma$, that is
\begin{equation}
\label{eq:periodN}
\frac{1}{2\pi\i} U=\int_{\infty^-}^{\infty^+}\omega=M_1+\mathbb B M_2,
\end{equation}
where $M_1, M_2\in\mathbb Q^g$ are column vectors, such that $NM_1, NM_2\in\mathbb Z^g$.
\end{theorem}
We can now use M\"obius transformations of the $z$-plane to transform equation \eqref{eq:Gamma} to equation \eqref{T}.
%
\begin{theorem}
\label{th:todaNdef}
Consider an $N$-periodic in $n$ solution  to the Toda lattice \eqref{eq:toda} constructed by formulas \eqref{eq:solTodav} and \eqref{eq:solTodac} from the hyperelliptic surface $\Gamma_{\x_0}$ of genus $g$ corresponding to the curve \eqref{T}, \eqref{delta} with $\x=\x_0$. For a value of $\x_0$ away from a set of measure zero defined in Theorem \ref{thm_converse-Toda}, there exists a continuous $g$-parameter deformation of this solution which remains $N$-periodic in $n$.
Moreover, any continuous deformation of this solution constructed from a family of curves $\Gamma_\x, \;\x\in\mathcal X,$ obtained from $\Gamma_{\x_0}$ by varying $\x=(x_1, \dots, x_g)$ and allowing  $\u=(u_1, \dots, u_g)$ to be functions of $\x$ and which remains $N$-periodic in $n$
solution to the Toda lattice is obtained by formulas \eqref{eq:solTodav} and \eqref{eq:solTodac} relative to the curves $\Gamma_\x, \;\x\in\mathcal X$. In this case, the  branch points $u_1, \dots, u_g$ of the coverings $u:\Gamma_\x\to\mathbb CP^1$ as functions of $\x$ have the derivatives expressed by \eqref{umder_g} and satisfy system \eqref{um_xkxn}, \eqref{um_xkxk} of Theorem \ref{thm_support}. The differentials $\Omega$ and $v_m$ in \eqref{umder_g} are defined on $\Gamma_\x$ by \eqref{Omega} and \eqref{v}, respectively.
\end{theorem}

{\it Proof.} If $(v_n, c_n)$ is a  periodic  in $n$ solution to the Toda lattice equation with period $N$ constructed by \eqref{eq:solTodav} and \eqref{eq:solTodac} on $\Gamma_{\x_0}$, then, according to Theorem \ref{th:todaN} the surface $\Gamma_{\x_0}$ satisfies Definition \ref{def_Toda} after relaxing the condition of reality for the branch points $x_j$ and $u_j$. Let $\Omega_0(\x_0)$ be the differential defined on $\Gamma_{\x_0}$ by \eqref{Omegalpha} with $\alpha=0$. Its vector of $b$-periods \eqref{b-Omegalpha} coincides with $U$ \eqref{eq:periodN}.
For a generic $\x_0$, Theorem \ref{thm_converse-Toda} shows the existence of a local variation of $\x=(x_1, \dots, x_g)$ entailing a deformation $\Gamma_\x,\;\x\in\mathcal X,$ of the surface $\Gamma_{\x_0}$  in such a way that $b$-periods of $\Omega_0(\x)$ stay constant for $\Omega_0(\x)$ being the differential \eqref{Omegalpha} with $\alpha=0$ on the surface $\Gamma_\x$. Thus condition \eqref{eq:periodN} remains satisfied for the surfaces $\Gamma_\x,\;\x\in\mathcal X$. This provides a continuous deformation of the solution $(v_n, c_n)$ given by \eqref{eq:solTodav}, \eqref{eq:solTodac} for the curves $\Gamma_\x$ preserving the period $N$ of $(v_n, c_n)$.

Let now  $(v_n(\x), c_n(x)),\;\x\in\mathcal X,$ be a continuous deformation of $(v_n, c_n)$ which remains $N$-periodic in $n$ and which is constructed from the family of surfaces $\Gamma_\x$ obtained from $\Gamma_{\x_0}$ by varying $\x$ and allowing  $\u$ to be functions of $\x$.  Theorem \ref{th:todaN} implies that for each $\x$, the solution $(v_n(\x), c_n(\x))$
 is obtained by formulas \eqref{eq:solTodav} and \eqref{eq:solTodac} with the underlying surface $\Gamma_\x$ satisfying condition \eqref{eq:periodN}.  We thus conclude that these underlying surfaces $\Gamma_\x, \;\x\in\mathcal X$ form a continuous Toda family in the sense of Definition \ref{def_Todafamily}. The rest of the proof follows from Theorems \ref{thm_umder} and  \ref{thm_support}. $\Box$
\\

Restricting to the case of real $\x, \u\in\mathbb R^g$, the surfaces $\Gamma_\x$ corresponding to the hyperelliptic curves \eqref{T}, \eqref{delta} are Toda curves and there is an integer $N$ and  a polynomial $\mathcal P_N$ satisfying the Pell equation  \eqref{Pell} corresponding to each $\Gamma_\x$.
Denoting $\mu^2=\Delta_{2g+2} \mathcal Q^2_{N-d}$ in \eqref{Pell}, we get that the curve $\Gamma_\x$ which produces $N$-periodic in $n$ solutions \eqref{eq:solTodav} and \eqref{eq:solTodac} of the Toda lattice \eqref{eq:toda} can be mapped to the curve of the following equation
\begin{equation}\label{eq:TodaNCurve}
\mathcal C_T: \mu^2=\mathcal P_N^2(z)-1,
\end{equation}
where $\mathcal P_N$ is a polynomial of degree $N$. The zeros $E_j$ of the polynomial $\mathcal P_N^2(z)-1$ are the end-points of the zones of the spectrum
of the operator $\mathbb {\tilde L}$, which is the restriction of $\mathbb L$ to the space of eigenfunctions of the operator of the shift by the period. Here $2\mathcal P_N$ is the trace of the monodromy matrix, see e.g. \cite{Nov}.  The polynomial $\mathcal P_N^2(z)-1$ has double roots exactly at the roots of the polynomial $\mathcal Q_{N-d}$. They form the singular, double points of the curve $\mathcal  C_T$,   each of them being an image of two points from $\Gamma_\x$.

Conversely, let us start from a curve $\mathcal C_T$ of the form \eqref{eq:TodaNCurve}, with a polynomial $\mathcal P_N$ of degree $N$ for which we assume  that $\mathcal P_N(z)^2-1$ has at most double zeros. Then, we can collect all the double zeros into a polynomial $\mathcal Q_{N-d}(z)$ of degree $N-d$ and obtain the relation
\begin{equation}\label{eq:TodaPell}
\mathcal P_N^2(z)-1=\Delta_{2d}(z)\mathcal Q_{N-d}^2(z),
\end{equation}
where $\Delta_{2d}$ is a polynomial of degree $2d$ with only simple zeros.   Introducing $v=\mu/\mathcal Q_{N-d}(z)$ we have a map from $\mathcal C_T$
to the Toda curve of genus $d-1$ given by the equation $v^2=\Delta_{2d}(z),$ which is essentially the curve \eqref{T}, \eqref{delta}.  Theorem \ref{th:todaNdef} also governs deformations of the curves \eqref{eq:TodaNCurve},  which, according to Proposition \ref{prop:isoequirat} and Corollary \ref{cor:isoequi}, produces in this case real-valued functions $u$ of real variables $\x$. We are going to exploit this further in Sections \ref{sect_kdvdifference} and \ref{sect_SW}.
 Observe that all the zeros of   $\Delta_{2d}$ and $Q_{N-d}$ are real if and only if  $\mathcal {P}_N$ is a generalized Chebyshev polynomial.

\

\section{The Korteweg - de Vries difference equation and Chebyshev polynomials}
\label{sect_kdvdifference}

The Korteweg - de Vries difference equation \cite{DMN} is an instance of the Toda lattice, with $v_n=0$ for all $n$.
The equation is given in the form
\begin{equation}\label{eq:kdvdifference}
\dot c_n=c_n(c_{n+1}-c_{n-1}), \qquad n\in \mathbb Z.
\end{equation}

According to \cite{DMN}, solutions to the KdV difference equation \eqref{eq:kdvdifference} are obtained  by \eqref{eq:solTodac} from  solutions to the Toda lattice  \eqref{eq:solTodac} for the Riemann surfaces $\Gamma$ of the form
\begin{equation}
\label{eq:kdvsurf}
\mu^2=\prod_{j=1}^{g+1}(z^2-z_i^2),\qquad 0<z_1<z_2<\dots<z_{g+1},
\end{equation}
i.e. such that the  set of branch points is real and symmetric with respect to $z=0$. Thus, a solution $c_n$ of the Korteweg - de Vries difference equation is $N$-periodic in $n$ if and only if the solution $(c_n, v_n=0)$ of the Toda lattice is $N$-periodic in $n$.
 From Theorem \ref{th:todaN}  and the considerations from the last part of Section \ref{sect_lattices}, in particular \eqref{eq:TodaPell} and reality of $\pm z_j$,   we get that the $N$-periodic solutions to the
KdV difference equation correspond  to supports of Pell's equation symmetric with respect to the origin,  supporting a generalized Chebyshev polynomial $\mathcal P_N$ of degree $N$.
In particular, we have $M_2=0$ in this case in \eqref{eq:periodN}, in the accordance with Corollary \ref{cor:M20}.
Starting form a periodic solution to the KdV difference equation \eqref{eq:kdvdifference} of period $N$, we want to study deformations of the curve \eqref{eq:kdvsurf} which would produce new periodic solution to the KdV difference equation of the same period. We cannot apply isoequilibrium deformations directly, because {\it a priori} we cannot guarantee that the deformed intervals would remain symmetric with respect to the origin. Thus, to that end, we need to provide a more detailed analysis of periodic solutions to the KdV difference equation, which we do next.

\begin{proposition}\label{prop:oddgen}
\begin{itemize}
\item [(i)] For  odd genera $g=2k-1$  of the curves \eqref{eq:kdvsurf}, the periodic solutions  to the KdV difference equation \eqref{eq:kdvdifference} are of even period $N=2\ell$ in $n$ and the generalized Chebyshev polynomials $\mathcal P_{2\ell}(z)$ corresponding to them are even. There are no periodic solutions with odd period in this case.

 \item [(ii)] If a generalized Chebyshev polynomial $\mathcal P_{2\ell}(z)$ is associated with $2\ell$ periodic solution  to the KdV difference equation constructed using the curve \eqref{eq:kdvsurf} of genus $g=2k-1$, then the polynomial $\hat p_{\ell}(w)$ of degree $\ell$ given by $\hat p_{\ell}(w)=\mathcal P_{2\ell}(z)$ with $w=z^2$, satisfies the Pell equation
     \begin{equation}\label{eq:pellkdv}
     \hat p_{\ell}^2(w) - \hat \Delta_{2k}(w)\hat q_{\ell-k}^2(w)=1,
     \end{equation}
     where
\begin{equation}
\label{delta_pellkdv}
\hat \Delta_{2k}(w)=\prod_{j=1}^{g+1=2k}(w-z_j^2)
\end{equation}
     is a polynomial of degree $2k$ with positive zeros  and $\hat q_{\ell-k}$ is a polynomial of degree $\ell-k$.

      \item [(iii)] If the polynomial $\hat {p}_{\ell}(w)$ of degree $\ell$ satisfies the Pell equation \eqref{eq:pellkdv}
     where $\hat \Delta_{2k}(w)$ \eqref{delta_pellkdv}
     is a polynomial of degree $2k$ with positive zeros, then  $\mathcal P_{2\ell}(z)=\hat {p}_{\ell}(w)$ with $w=z^2$ is  the generalized Chebyshev polynomial associated with the $2\ell$ periodic solution to the KdV difference equation constructed using the curve \eqref{eq:kdvsurf} of genus $g=2k-1$.
     \end{itemize}
\end{proposition}

{\it Proof.}
\begin{itemize}
\item[(i)]  Let the curve \eqref{eq:kdvsurf} correspond to a support of the Pell equation.  Since the support is symmetric with respect to the origin, if $\mathcal P_N(z)$ is the associated generalized Chebyshev polynomial, so is $(-1)^N\mathcal P_N(-z)$. Thus, from the uniqueness of the associated Chebyshev polynomial, it follows that $\mathcal P_N(z)$ is odd for $N$ odd and it is even for $N$ even. Thus, for $N$ odd $\mathcal P_N(z)$ is odd and $\mathcal P_N(0)=0$. For odd genera, the interval containing $0$ is a gap interval, and this leads to the contradiction.

\item[(ii), (iii)] Suppose now that $N=2\ell$ is even and a polynomial $\mathcal P_{2\ell}(z)$ is even.  Then there  exists a polynomial $\hat p_{\ell}$ of degree $\ell$, such that
$\mathcal P_{2\ell}(z)=\hat p_{\ell}(w)$, with $w=z^2$. One can observe that $\mathcal P^2_{2\ell}(z)-1$ has $4\ell$ real zeros if and only if $\hat p^2_{\ell}(w)-1$ has $2\ell$ real  positive zeros. We get that
$\hat p_{\ell}(w)$ satisfies \eqref{eq:pellkdv},
     where $\hat \Delta_{2k}$ \eqref{delta_pellkdv}
    is a polynomial of degree $2k$ with positive zeros, if and only if
$$
\mathcal {P}_{2\ell}^2(z) - \prod_{j=1}^{g+1=2k}(z^2-z_j^2)\mathcal {Q}_{2\ell-2k}^2(z)=1,
$$
for some polynomial $\mathcal {Q}_{2\ell-2k}$ of degree $2\ell-2k.$
One should observe that in this case $\mathcal {Q}_{2\ell-2k}(z)$ is even: given that the square of $\mathcal {Q}_{2\ell-2k}(z)$ is even, which follows from the Pell equation, the polynomial $\mathcal {Q}_{2\ell-2k}(z)$ is either even or odd. But since its degree is even, it cannot be odd.
\end{itemize}
$\Box$

\begin{proposition}\label{prop:evengen}
\begin{itemize}
\item [(i)]
For even genera $g=2k$  of the curves \eqref{eq:kdvsurf}, if the periodic solutions to the KdV difference equation \eqref{eq:kdvdifference} are of odd period $N=2\ell+1$ then the generalized Chebyshev polynomials $\mathcal P_{2\ell+1}$ corresponding to them  are odd; if the periodic solutions to the KdV difference equation \eqref{eq:kdvdifference} are of even period $N=2\ell$ then the generalized Chebyshev polynomials $\mathcal P_{2\ell}$ corresponding to them are even.
\item [(ii)]
For even genera $g=2k$ and odd period $N=2\ell+1$, there exists a polynomial $\hat p_{\ell}(w)$ of degree $\ell$ such that $ \mathcal P_{2\ell+1}(z)=z\hat p_{\ell}(w)$, with $w=z^2$. The polynomial $\hat p_{\ell}(w)$ satisfies a modified Pell's equation
\begin{equation}\label{eq:Pelmod}
w\hat p_{\ell}^2(w)-\prod_{j=1}^{2k+1}(w-z_j^2)\hat q_{\ell-k}^2(w)=1,
\end{equation}
 where $\hat q_{\ell-k}$ is a polynomial of degree $\ell-k.$
The polynomial
\begin{equation}\label{eq:tildep}
\tilde p_{2\ell+1}(w)=2w\hat p^2_{\ell}(w)-1
\end{equation}
 satisfies Pell's equation:
\begin{equation}\label{eq:PelPel}
\tilde p_{2\ell+1}^2(w)-4w\prod_{j=1}^{2k+1}(w-z_j^2)\hat p_{\ell}^2(w)\hat q_{\ell-k}^2(w)=1.
\end{equation}

\item [(iii)]
For the even genera $g=2k$  of the curves \eqref{eq:kdvsurf} and even period $N=2\ell$,  there exists a polynomial $\hat p_{\ell}(w)$ of degree $\ell$ such that $ \mathcal P_{2\ell}(z)=\hat p_{\ell}(w)$, with $w=z^2$.  The polynomial $\hat p_{\ell}(w)$ satisfies Pell's equation
\begin{equation}\label{eq:Peleven}
\hat p_{\ell}^2(w)-w\prod_{j=1}^{2k+1}(w-z_j^2)\hat q_{\ell-k-1}^2(w)=1,
\end{equation}
where $\hat q_{\ell-k-1}$ is a polynomial of degree $\ell-k-1$ if and only if  $\mathcal P_{2\ell}(z)$ is an even generalized Chebyshev polynomial.
\end{itemize}
\end{proposition}

{\it Proof.} The first part follows along the same lines as the initial part of the proof of Proposition \ref{prop:oddgen}.
\begin{itemize}
\item[(ii)] Let now $g=2k$ and $N=2\ell+1$. Since $\mathcal P_{2\ell+1}(z)$ is odd, there exists a polynomial $\hat p_{\ell}(w)$ of degree $\ell$ such that
$ \mathcal P_{2\ell+1}(z)= z\hat p_{\ell}(w)$, with $w=z^2$.

From Pell's equation for $\mathcal P_{2\ell+1}(z)$,
$$
\mathcal P_{2\ell+1}^2(z) - \prod_{j=1}^{2k+1}(z^2-z_j^2)\mathcal Q_{2(\ell-k)}^2(z)=1,
$$
with the observation that $ \mathcal Q_{2(\ell-k)}(z) $ is an even function (using an argument similar to the one at the end of the proof of Proposition \ref{prop:oddgen}) and defining $q_{\ell-k}(w)=\mathcal Q_{2(\ell-k)}(z)$, we get
$$
w\hat p_{\ell}^2(w)-\prod_{j=1}^{2k+1}(w-z_j^2)\hat q_{\ell-k}^2(w)=1.
$$
One verifies directly that the polynomial  $\tilde p_{2\ell+1}(w)$, defined by \eqref{eq:tildep},
 satisfies Pell's equation \eqref{eq:PelPel}.
\item[(iii)] Let $g=2k$ and $N=2\ell$. Since $\mathcal P_{2\ell}(z)$ is even, there exists a polynomial $\hat p_{\ell}(w)$ of degree $\ell$ such that
$ \mathcal P_{2\ell}(z)= \hat p_{\ell}(w)$, with $w=z^2$. From Pell's equation for $\mathcal P_{2\ell}(z)$:
$$
\mathcal P_{2\ell}^2(z) - \prod_{j=1}^{2k+1}(z^2-z_j^2)\mathcal Q_{2(\ell-k)-1}^2(z)=1
$$
with the observation that $ \mathcal Q_{2(\ell-k)-1}(z) $ is an odd function (using an adjustment to the argument  at the end of the proof of Proposition \ref{prop:oddgen})  and defining $z\hat q_{\ell-k-1}(w)=\mathcal Q_{2(\ell-k)-1}(z)$, we get  \eqref{eq:Peleven}.  The converse follows directly.
\end{itemize}
$\Box$

Now, we are ready to study the main question of this section: {\it given a periodic solution to the KdV difference equation \eqref{eq:kdvdifference} of period $N$, how to deform the  underlying curve \eqref{eq:kdvsurf}  in order for the deformed curves to give rise to  new periodic solution to the KdV difference equation of the same period.}

\begin{theorem}\label{thm:kdvdef}
Given a periodic solution to the KdV difference equation \eqref{eq:kdvdifference} of period $N$,  the following deformations of the underlying hyperelliptic curve \eqref{eq:kdvsurf}  produce solutions  $c_n$, given by \eqref{eq:solTodac}, to the same KdV difference equation of the same period:

\begin{itemize}
\item[(i)] For $g=2k-1$, $N=2\ell$, we consider isoequilibrium deformations  of the set of non-negative intervals obtained from the set of positive intervals
\begin{equation}
\label{positiveintervals}
[z_1^2, z_2^2], \; [z_3^2, z_4^2], \dots, [z_{2k-1}^2, z_{2k}^2]
\end{equation}
by an affine transformation that maps $z^2_1, z^2_2$ to $0, 1$.
This corresponds to isoperiodic deformations of Toda curves of genus $g_1=k-1$,  defined as hyperelliptic coverings $\Tx$ of $\mathbb CP^1$ having the intervals \eqref{positiveintervals} as branch cuts. In this case,  the  branch points $u_1=z_4^2, \dots, u_{k-1}=z_{2k}^2$ of the coverings $\Tx$, as functions of $\x=(z_3^2, z_5^2, \dots, z_{2k-1}^2)$, have derivatives expressed by \eqref{umder_alpha}, \eqref{v} with $\alpha=0$
and satisfy system \eqref{um_xkxn}, \eqref{um_xkxk} of Theorem \ref{thm_support}, after an affine transformation that maps $z^2_1, z^2_2$ to $0, 1$, respectively.
The equilibrium deformations of the set of intervals \eqref{positiveintervals}  naturally generate continuous deformations of the intervals
\begin{equation}
\label{intervals_th8}
[-z_{2k}, - z_{2k-1}], \dots, [-z_4,-z_3], [-z_2, -z_1], [z_1, z_2], [z_3,z_4], \dots, [z_{2k-1}, z_{2k}],
\end{equation}
throughout which the intervals \eqref{intervals_th8} remain the support of a generalized Chebyshev polynomial of degree $2\ell$.

\item[(ii)] For $g=2k$, $N=2\ell$, we consider isoequilibrium deformations of the set of  non-negative intervals obtained from the set of  intervals
\begin{equation}
\label{nonnegativeintervals1}
[0, z_1^2], [z_2^2, z_3^2], \dots, [z_{2k}^2, z_{2k+1}^2],
\end{equation}
by an affine transformation that maps $0, z^2_1$ to $0, 1$. These deformations  keep $0$ as the left end-point of the leftmost interval in \eqref{nonnegativeintervals1}. They also generate naturally deformations of the intervals
\begin{equation}
\label{intervals_th8ii}
[-z_{2k+1}, - z_{2k}], \dots, [-z_3,-z_2], [-z_1, z_1], [z_2,z_3], \dots, [z_{2k}, z_{2k+1}],
\end{equation}
throughout which the intervals \eqref{intervals_th8ii} remain the support of a generalized Chebyshev polynomial of degree $2\ell$.

\item[(iii)] For $g=2k$, $N=2\ell+1$, we consider isoequilibrium deformations of the set of  non-negative intervals obtained from the set of  intervals
\begin{equation}
\label{nonnegativeintervals2}
[0, z_1^2], [z_2^2, z_3^2], \dots, [z_{2k}^2, z_{2k+1}^2],
\end{equation}
by an affine transformation that maps $0, z^2_1$ to $0, 1$. These deformations keep $0$ as the left end-point of the leftmost interval in \eqref{nonnegativeintervals2} and generate naturally deformations of the intervals
\begin{equation}
\label{intervals_th8iii}
[-z_{2k+1}, - z_{2k}], \dots, [-z_3,-z_2], [-z_1, z_1], [z_2,z_3], \dots, [z_{2k}, z_{2k+1}],
\end{equation}
throughout which the intervals \eqref{intervals_th8iii} remain the support of a generalized Chebyshev polynomial of degree $2\ell+1$.
\end{itemize}

For (ii) and (iii), this corresponds to isoperiodic deformations of Toda curves of genus $g_1=k$  defined as hyperelliptic coverings $\Tx$ of $\mathbb CP^1$ having the intervals \eqref{nonnegativeintervals1} or \eqref{nonnegativeintervals2}, respectively, as branch cuts. In these cases,  the  branch points $u_1=z_3^2, \dots, u_{k}=z_{2k+1}^2$ of the coverings $\Tx$ as functions of $\x=(z_2^2, z_4^2, \dots, z_{2k}^2)$ have derivatives expressed by \eqref{umder_alpha}, \eqref{v} with $\alpha=0$
and satisfy system \eqref{um_xkxn}, \eqref{um_xkxk} of Theorem \ref{thm_support}, after an affine transformation that maps $0, z^2_1$ to $0, 1$, respectively.
\end{theorem}

{\it Proof.}  Note that isoequlibrium deformations of a set of real intervals exist due to Proposition \ref{prop:isoequiirat}. Further proof  of (i) follows from Proposition \ref{prop:oddgen},  Theorem \ref{cor:rational}, Proposition \ref{prop:isoequirat},  and Corollary \ref{cor:isoequi}.

Proofs of (ii) and (iii) follow from Proposition \ref{prop:evengen}, Theorem \ref{cor:rational}, Proposition \ref{prop:isoequirat}, Corollary \ref{cor:isoequi}, and,
in the case (iii), the following argument.

With the initial curve \eqref{eq:kdvsurf} is associated a generalized Chebyshev polynomial $\mathcal P_{2\ell+1}$ and a polynomial $\hat p_{\ell}$ defined by $ \mathcal P_{2\ell+1}(z)= z\hat p_{\ell}(w)$, with $w=z^2$. The deformations of the set of intervals \eqref{nonnegativeintervals2} is such that the polynomial $\tilde p_{2\ell+1}$ defined by \eqref{eq:tildep} for the initial set of intervals admits a corresponding continuous deformation throughout which it satisfies the Pell equation \eqref{eq:PelPel} for the deformed set of $z_j^2$ for some polynomial $\mathcal Q_{2\ell-k}(w)$ replacing $p_{\ell}(w)\hat q_{\ell-k}(w)$. From \eqref{eq:PelPel}, we have the following factorization for the initial situation:
\begin{equation}
\label{factorization}
\tilde p_{2\ell+1}^2(w)-1=(\tilde p_{2\ell+1}(w)+1)(\tilde p_{2\ell+1}(w)-1) = 4w\prod_{j=1}^{2k+1}(w-z_j^2)\hat p^2_{\ell}(w)\hat q_{\ell-k}^2(w).
\end{equation}
From continuity of our deformations, from the continuous dependance of the polynomials on their zeros and vice-versa, and from \eqref{factorization}, we conclude that the polynomial  $\hat p_{\ell}(w)$ is unambiguously defined for a deformed set of intervals, for small enough deformations. Indeed, the zeros of the polynomial $\hat p_{\ell}(w)\hat q_{\ell-k}(w)$ are simple, which means that the zeros of $\hat p_{\ell}(w)$ are initially simple and distinct from the zeros of $\hat q_{\ell-k}(w)$ (as well as from the endpoints of the deformed intervals from the set \eqref{nonnegativeintervals2}). We deform the intervals in a way that they remain being the support of a generalized Chebyshev polynomial. That means that the polynomial $\hat q_{\ell-k}(w)$ continuously deforms, inducing continuous deformations of its zeros. As long as the deformations are small enough that the deformed zeros remain simple, the deformed polynomial  $\hat p_{\ell}(w)$ is well defined.

 Thus, under small continuous deformations, $\tilde p_{2\ell+1}(w)+1$ keeps being the product of  the same (up to the deformation) factors of the right hand side of \eqref{factorization}. Therefore, since \eqref{eq:tildep} is satisfied for the initial situation, the same relation, up to possible small deformation of the  scalar factor $\delta>0$
\begin{equation}
\label{pp}
\tilde p_{2\ell+1}(w)+1=\delta w\hat p^2_{\ell}(w)
\end{equation}
holds after small deformations.

Now, from relation \eqref{pp} and
$$
\tilde p_{2\ell+1}(w)-1=\frac{4}{\delta} \prod_{j=1}^{2k+1}(w-z_j^2)\hat q_{\ell-k}^2(w),
$$
we get
\begin{equation}
\label{Pelldeformation}
\frac{\delta}{2}w\hat p^2_{\ell}(w)-\frac{2}{\delta}\prod_{j=1}^{2k+1}(w-z_j^2)\hat q_{\ell-k}^2(w)=1.
\end{equation}
This relation, up to nonessential factors, which can be absorbed in the polynomials, coincides with  \eqref{eq:Pelmod}.
We can now define $ \mathcal P_{2\ell+1}(z)= \sqrt{\delta/2}z\hat p_{\ell}(w)$, with $w=z^2$. Equation \eqref{Pelldeformation} becomes the Pell equation for $ \mathcal P_{2\ell+1}(z)$ and thus, $ \mathcal P_{2\ell+1}(z)$  is the generalized Chebyshev polynomial of degree $2\ell+1$ supported by the set of intervals from \eqref{nonnegativeintervals2}. This
 finishes the proof.
$\Box$

\section{Applications to the $SU(N)$ Seiberg-Witten theory}
\label{sect_SW}

In their foundational works \cite{SW1, SW2}, Seiberg and Witten considered the construction of the non-perturbative dynamics of $\mathcal N=2$ $SU(2)$-invariant super-symmetric Yang-Mills theory in the limit of low energy and momenta. The key ingredient in their construction was a family of elliptic curves and meromorphic differentials on them, along with possible degenerations of elliptic curves into rational curves with singularities.

Hyperelliptic curves of the form \eqref{eq:TodaNCurve} of genus $N-1$ play a crucial role in the $SU(N)$ gauge group generalizations of the Seiberg-Witten theory, as discovered in \cite{AF1995, KLTY}. It was observed in \cite{GKMMM} that the curves introduced in \cite{AF1995} and \cite{KLTY} coincide with the curves that had appeared before in the studies of the Toda lattices. A deeper connection of the Seiberg-Witten theory with the so-called Witham-Toda hierarchy was developed in \cite{NaTa}.

In a special case of the $\mathcal N=2$ $SU(N)$ super-symmetric Yang-Mills theory, the so-called case without fundamental hypermultiplets (see e.g. \cite{DHPh}),
the main object is the family of curves
\begin{equation}
\label{eq:SWCurve}
\mathcal C_T(x_1,\dots, x_n): \mu^2=\mathcal P_N^2(z)-\hat\Lambda^2,
\end{equation}
 parametrized by $n$ complex parameters $x_1,\dots, x_n$,  where $\hat\Lambda$ is a real constant. Here $\mathcal P_N(z)$ is a polynomial of degree $N$. The parameters play the role of vacuum moduli  of the Yang-Mills theory. They can be chosen as a subset
of the set of zeros of the varying  polynomial  $\mathcal P_N^2(z)-\hat \Lambda^2$. In general, these curves are non-singular  hyperelliptic curves of genus
$n=N-1$, which corresponds to the case when all zeros of $\mathcal P_N^2(z)-\hat\Lambda^2$ are simple. However, for certain values of the parameters, the curves may become singular,  as some of the zeros of $\mathcal P_N^2(z)-\hat\Lambda^2$ merge to form a double zero. When such singularities occur, some of the particles in the Yang-Mills theory become massless. We will call such situations {\it singular regimes}   and denote by $N-g-1$ the number of double zeros of $\mathcal P_N^2(z)-\hat\Lambda^2$.  We call such particles {\it new massless particles.}

In other words, in a singular regime, coming with $N-g-1$ double zeros of $\mathcal P_N^2(z)-\hat\Lambda^2$, there exists a polynomial $\Delta_{2g+2}(z)$ of degree $2g+2$ with only simple zeros  (coinciding with the simple zeros of $\mathcal P_N^2(z)-\hat\Lambda^2$) and a polynomial $\mathcal Q_{N-g-1}(z)$ of degree $N-g-1$ with only simple zeros, such that the polynomial $\mathcal P_N$ from \eqref{eq:SWCurve} satisfies Pell's equation \eqref{Pell} (assuming that we made a nonessential modification of the Pell's equation by  replacing 1 in the right hand side by the constant $\hat\Lambda^2$). Thus, in the case when all zeros of $\Delta_{2g+2}(z)$  and $\mathcal Q_{N-g-1}(z)$ are real,  $\mathcal P_N$ is a generalized Chebyshev polynomial.  We will also consider in this section the situations when some of the polynomials  $\Delta_{2g+2}(z)$  and $\mathcal Q_{N-g-1}(z)$ have some of the zeros which are not real.

The desingularization of the singular curve \eqref{eq:SWCurve} is a hyperelliptic curve $\Gamma_\x$ of genus $g$ defined by the equation
\begin{equation}
\label{eq:SWdesing}
\Gamma_\x: w^2=\Delta_{2g+2}(z).
 \end{equation}
where $\x=(x_1, \dots, x_g)$ is a subset of the set of zeros of $\Delta_{2g+2}.$
New massless particles correspond to $N-g-1$  zeros of the polynomial $\mathcal Q_{N-g-1}$. In the case when all zeros of $\Delta_{2g+2}(z)$  and $\mathcal Q_{N-g-1}(z)$ are real,  the zeros of the polynomial $\mathcal Q_{N-g-1}$ are the {\it internal} critical points of the generalized Chebyshev polynomial $\mathcal P_{N}$, that is  the critical points lying inside the intervals of the support of $\mathcal P_{N}$ (see the proof of Theorem \ref{cor:rational}).

\begin{example}\label{ex:ellipticmassless} Let us consider one example from the  $SU(2)$ case.
Seiberg and Witten in \cite{SW2} considered the  curve of the equation
\begin{equation*}
 Y^2=X^3+2uX^2+\Lambda^4X,
\end{equation*}
where $\Lambda$ is a real constant.
By a M\"obius transformations in $X$, this curve can be brought to the form \cite{AF1995}:
\begin{equation*}
\mathcal C_T(u): \mu^2=\Big(z^2-\frac{1}{2}u\Big)^2-\Lambda^4.
\end{equation*}

This is an elliptic curve for $u\ne \pm 2\Lambda^2$.
The polynomial $\mathcal P_2$ in this example is
$$\mathcal P_2(z)=z^2-\frac{u}{2}.$$

Let us consider the singular case $u=2\Lambda^2$. Then the equation of the the curve becomes
$\mu^2=(z^2-\Lambda^2)^2-\Lambda^4,$
which is a singular cubic
$$
\mathcal  C_T(2\Lambda^2): \mu^2=z^2(z^2-2\Lambda^2).
$$
With $(z, \mu)$ mapping to $(z, w)$ by $w=\mu/z$, we get a nonsingular rational curve:
$$
\Gamma: w^2=z^2-2\Lambda^2.
$$
The polynomial $\mathcal P_2$ in this singular case is thus
\begin{equation*}
\mathcal P_2(z)=z^2-\Lambda^2.
\end{equation*}
This is the monic  Chebyshev polynomial of degree two corresponding to the interval $[-\Lambda \sqrt {2}, \Lambda \sqrt {2}]$.
It satisfies Pell's equation
\begin{equation*}
\mathcal P_2(z)^2-w^2z^2=\Lambda^4,
\end{equation*}
where $\Delta_{2}(z)=w^2=z^2-2\Lambda^2$ and $\mathcal Q_{1}(z)=z$. The internal critical point of $\mathcal P_2$, which is $z=0$, the zero of polynomial $\mathcal Q_1$, corresponds to one massless particle which arises in this singular case. For $\Lambda = \sqrt{2}/2$ we get  $\mathcal P_2(z)=z^2-1/2$ which is the classical monic Chebyshev polynomial of degree two on the interval $[-1, 1]$.

\end{example}

In what follows, we refer to the parameters $\x=(x_1, \dots, x_g)$ from \eqref{eq:SWdesing} as the {\it vacuum moduli parameters} of the singular regime.
\begin{theorem}\label{th:singulardef}
Consider the vacuum moduli parameters $\x_0=(x^0_1, \dots, x^0_g)$  of a singular regime of the Yang-Mills theory  with $N-g-1$ double zeros of $\mathcal P_N^2-1$. Let $\Delta^0_{2g+2}$ of \eqref{eq:SWdesing} be given by \eqref{delta} with distinct $x^0_j, u^0_j\in\mathbb C.$ For generic moduli parameters $\x_0$ lying away from some set of measure zero defined in Theorem \ref{thm_converse-Toda}, there exists a local continuous deformation of this singular regime which fixes the number of double zeros of $\mathcal P_N^2-1$.
This deformation is constructed from a family of curves $\Gamma_\x$ given by \eqref{eq:SWdesing} with $\x$ varying in some neighbourhood of $\x_0$ where  the  branch points $u_1, \dots, u_g$ of the coverings $z:\Gamma_\x\to\mathbb CP^1$ as functions of $\x=(x_1, \dots, x_g)$ have derivatives expressed by
\eqref{umder_g}
and satisfy system \eqref{um_xkxn}, \eqref{um_xkxk} of Theorem \ref{thm_support}. The differentials $\Omega$ and $v_m$ in \eqref{umder_g} are defined on $\Gamma_\x$ by \eqref{Omega} and \eqref{v}, respectively.

\end{theorem}
{\it Proof.} The above discussion shows that a singular regime of the Yang-Mills theory  corresponds to a polynomial $\mathcal P_N$  for which the polynomial $\mathcal P_N^2-1$ has $N-g-1$ double zeros and $2g+2$ simple zeros. In other words, we have the polynomial relation $\mathcal P_N^2-1 = \Delta^0_{2g+2}\mathcal Q^2$ of the form of the Pell equation. This implies, as explained at the end of Section \ref{sect_TodaPell}, that the curve $\Gamma_{\x_0}$ defined by \eqref{eq:SWdesing} with $\Delta_{2g+2}=\Delta^0_{2g+2}$ is a Toda curve of order $N$ in the sense of Definition \ref{def_Toda} if we relax the requirement of reality of $x^0_j$ and $u^0_j$. Therefore, \eqref{Abel} holds on $\Gamma_{\x_0}$ and one can define the differential $\Omega=\Omega(\x_0)$ by \eqref{Omega}.
Theorem \ref{thm_converse-Toda} shows the existence of a local continuous generalized Toda family $(\Gamma_{\x}, \Omega(\x))$ providing an isoperiodic deformation of the pair $(\Gamma_{\x_0}, \Omega(\x_0))$ by varying the parameters $\x=(x_1, \dots, x_g)$ in a neighbourhood of ${\x_0}$. Thus $\Gamma_{\x}$ is the curve defined by \eqref{eq:SWdesing} with $\Delta_{2g+2}$ obtained from $\Delta^0_{2g+2}$ by replacing $\x_0$ by $\x.$  The periods of $\Omega(\x)$ \eqref{Omega}  are given by $M_1$ and $M_2$ from \eqref{Abel}, see \eqref{Omega} and \eqref{bOmega}.
Since the periods of $\Omega(\x)$ are independent of $\x$, we have that \eqref{Abel} is satisfied for the surfaces $\Gamma_\x$ with constant rational vectors $M_1, M_2$ such that $NM_1, NM_2\in\mathbb Z^g.$ By the Abel theorem, this implies that $\Gamma_\x$ is a Toda curve of order $N$ for all $\x$ in the neighbourhood in question where we relax the requirement of the reality of the branch points.
  Therefore, $\mathcal P_N^2-1 = \Delta_{2g+2}\mathcal Q^2$ is satisfied for all $\x$ in some neighbourhood of $\x_0$. Here the polynomials $\mathcal P_N, \;\Delta_{2g+2}, \;\mathcal Q$ vary smoothly as $\x$ varies. This implies that for all such $\x$ the polynomial $\mathcal P_N^2-1$  has $2g+2$ simple zeros and $N-g-1$ double zeros, thus corresponding to  a singular regime of the Yang-Mills theory, for all $\x$ in some neighbourhood of $\x_0.$
This shows the existence of a local continuous deformation of the singular regime corresponding to the moduli parameters $\x_0$.
The rest of the statement  follows from Theorem \ref{thm_isoperiodicity}.
$\Box$
\\

 Under the assumption that the number of new massless particles corresponding to a given type of singular regime of the theory does not depend on small deformation of the vacuum moduli parameters $(x_1,\dots, x_n)$, it would follow from Theorem \ref{th:singulardef} that the number of new massless particles would remain unchanged under the above isoperiodic deformations.

In Example \ref{ex:ellipticmassless}, the set of the vacuum moduli parameters is empty ($g=0$) and there are no nontrivial deformations in this case.

In a more general case of the $\mathcal N=2$ $SU(N)$ super-symmetric Yang-Mills theory {\it with} fundamental hypermultiplets (see e.g. \cite{DHPh}),
the family of the associated curves is given by
\begin{equation*}
\mathcal C_T(x_1,\dots, x_n): \mu^2=\mathcal P_N^2(z)-\hat\Lambda^2R(z),
\end{equation*}
where $R(z)=\prod_{j=1}^{N_f}(z+m_j)$. Here $N_f$ is the number of multiplets and their masses are $m_j$, $j=1,\dots, N_f$.
We will consider the situation with multiplets in a separate publication.

\section{Triangular solution to constrained Schelsinger system.}
\label{sect_solution}

In this section, we give a solution  to the  constrained Schelsinger system  introduced in \cite{isoharmonic}. This is a restriction of the classical Schlesinger system to a submanifold in the space of independent variables. More precisely, for some integer $1\leq K\leq 2g-1$ denote $\B_K=\{0,1,x_1, \dots, x_{K}, u_1, \dots, u_{2g-K}\}$ so that for the set \eqref{branch} we have $\B=\B_g.$ Then for $a_j$ being an arbitrary element of the set  $\B=\{0,1,x_1, \dots, x_{K}, u_1, \dots, u_{2g-K}\}$,  the constrained Schlesinger system for matrices $\Aaj$ has the following form for $i=1, \dots, K$:
\begingroup
\allowdisplaybreaks
\begin{eqnarray}
\label{constrained}
&&\partial_{x_i} A_{a_j} = \frac{[A_{x_i}, A_{a_j}]}{x_i-a_j} + \sum_{k=1}^{2g-K}\frac{[A_{u_k}, A_{a_j}]}{u_k-a_j} \frac{\partial u_k}{\partial x_i}\,, \quad\mbox{for}\quad a_j\notin\{x_i, u_1, \dots, u_{2g-K}\};
\nonumber\\
&&\partial_{x_i} A_{u_m}\!\! = \frac{[A_{x_i}, A_{u_m}]}{x_i-u_m} +\!\!\!\!\sum_{\substack{k=1\\k\ne m}}^{2g-K} \!\frac{[ A_{u_k},A_{u_m}]}{u_k-u_m} \frac{\partial u_k}{\partial x_i}-\frac{\partial u_m}{\partial x_i}\!\!\sum_{\substack{a_j\in \B\\a_j\ne u_m}}\!\! \frac{[A_{a_j},A_{u_m}]}{a_j-u_m}\,,\;\;\mbox{for}\;\; 1\leq m\leq 2g-K;\qquad
\\
&&\partial_{x_i} A_{x_i} = -\sum_{\substack{a_j\in \B \\a_j\ne x_i}}\frac{[A_{x_i}, A_{a_j}]}{x_i-a_j} + \sum_{k=1}^{2g-K}\frac{[A_{u_k}, A_{x_i}]}{u_k-x_i} \frac{\partial u_k}{\partial x_i}\,.
\nonumber
\end{eqnarray}
\endgroup
The last equation may be replaced by $\sum_{a_j\in \B} A_{a_j}=const.$ If $K=2g$, equations \eqref{constrained} turn into the usual Schlesinger system.   The submanifold in the space of the variables of the classical Schlesinger system is given by setting the first two variables to $0$ and $1$ and by giving $u_1, \dots, u_{2g-K}$ as functions of $x_1, \dots, x_K.$ If, additionally, we have equations for the partial derivatives of these functions, then system \eqref{constrained} is well determined, otherwise it is underdetermined.

Let us set $K=g$ and consider system \eqref{constrained} together with equations \eqref{um_xkxn} and \eqref{um_xkxk} for the functions $u_j$. Such a system is well determined and a solution to it is given in  the next theorem  in terms of the meromorphic differentials defined over families of curves considered in this paper, where $\Aaj$ are $2\times 2$ upper triangular traceless matrices.

\begin{theorem}
\label{thm_css}
Let $\Tx\to X\subset \mathbb C^g$ be a family of Toda type as in Definition \ref{def_Todatype} carrying a meromorphic differential of the third kind $\Omega:=\Omega_\alpha$ defined by \eqref{Omegalpha} with some constant $\alpha\in\mathbb C^g$. Let the family of Riemann surfaces $\Tx$ correspond to the hyperelliptic curves defined by \eqref{surf} and $\x=(x_1, \dots, x_g)$ be a subset of branch points of the coverings $u:\Tx\to\mathbb CP^1$ as in \eqref{surf}.  Let  $\phi$ be the holomorphic differential \eqref{phi} defined on the surfaces from the family. Then, together with the branch points $u_1, \dots, u_g$ of the coverings $u$,  the following matrices satisfy the constrained Schlesinger system \eqref{constrained} with $K=g$ and system \eqref{um_xkxn}, \eqref{um_xkxk} with respect to the variables $x_1, \dots, x_g$.   Here  $a_j$ runs through the set $\B=\{0,1,x_1, u_1, \dots, x_g, u_g\}$  of all branch points of the coverings $u$:
\begin{equation}
\label{Aaj-mat}
A_{a_j}=\left(\begin{array}{cc}-\frac{1}{4} & \Aaj^{12} \\ 0 & \frac{1}{4}\end{array}\right)
\end{equation}
where
\begin{equation}
\label{Aaj}
A^{12}_{a_j} =  \frac{t}{4}\Omega(P_{a_j})\varphi(P_{a_j})
\end{equation}
with $t\in\mathbb C$ being an arbitrary constant.

\end{theorem}

{\it Proof.} Due to Theorem \ref{thm_isoperiodicity}, the branch points $u_1, \dots, u_g$ as functions of $\x$ satisfy system \eqref{um_xkxn}, \eqref{um_xkxk} and also \eqref{umder_g}. Using notation \eqref{Aaj}, we rewrite \eqref{umder_g} in the form
\begin{equation}
\label{umder_A}
 \frac{\partial u_m}{\partial x_i}=-\frac{\Axi^{12}}{\Aum^{12}} \prod_{\alpha=1, \alpha\neq m}^g\frac{x_i-u_\alpha}{u_m-u_\alpha}\,.
\end{equation}

The  arbitrary parameter $t$ in the solution in Theorem \ref{thm_css} reflects the invariance of the constrained Schlesinger system \eqref{constrained} under the following simultaneous, for all $a_j\in \B=\B_g,$ rescaling of the anti-diagonal elements of the matrices: $\Aaj^{12}\mapsto t\Aaj^{12} $ and $\Aaj^{21}\mapsto \frac{1}{t}\Aaj^{21}$.

T show that the sum of all the matrices $\Aaj$ \eqref{Aaj-mat} is constant, we only need to show that the sum of $\Aaj^{12}$ over all $a_j\in\B$ is a constant. This sum is actually zero as a sum of residues of the following differential on the compact Riemann surface $\Tx$ for some fixed $\x:$
\begin{equation}
\label{A12_g}
A^{12}(u)du := \frac{t}{2}\frac{\Omega(P)du(P)}{\varphi(P)\,v^2}.
\end{equation}
This differential has simple poles at $\paj$ with $a_j\in\B$ and no other singularities on $\Tx$. Moreover
 $A^{12}(u)$ is a well-defined function of $u\in\mathbb CP^1$ with simple poles in the set $B=\{0,1, x_1, \dots, x_g, u_1, \dots, u_{g}\}\subset\mathbb C$ and no other singularities. Thus, we can write
\begin{equation}
\label{A12_aj_g}
A^{12}_{a_j} =\underset{u=a_j}{\rm res} A^{12}(u) du =  \frac{t}{4}\Omega(P_{a_j})\varphi(P_{a_j})
\end{equation}
for any $a_j\in \B$, where the evaluation of differentials at ramification points is understood in the sense of  \eqref{evaluation}. Therefore, we have $\sum_{a_j\in B}A^{12}_{a_j} =0$ and
\begin{equation}
\label{sum_A12}
\sum_{a_j\in B}\Aaj =\left(\begin{array}{cc}-\frac{g+1}{2} & 0 \\ 0 & \frac{g+1}{2}\end{array}\right).
\end{equation}

It remains to prove that the matrices \eqref{Aaj-mat}, \eqref{Aaj} satisfy the first two equations of \eqref{constrained}.
For the matrices  $\Aaj$ from Theorem \ref{thm_css}, these equations reduce to the following
\begin{eqnarray}
\label{Aaj12}
&&\frac{\partial A_{a_j}^{12}}{\partial x_i} = \frac{1}{2}\frac{  A^{12}_{x_i}-A_{a_j}^{12} }{x_i-a_j} + \frac{1}{2} \sum_{k=1}^{g} \frac{ A^{12}_{u_k}- A_{a_j}^{12} }{u_k-a_j} \frac{\partial u_k}{\partial x_i}\,,
\qquad\qquad a_j\notin\{u_1, \dots, u_{g}\};\qquad\qquad\qquad
\\
&&\frac{\partial A_{u_m}^{12}}{\partial x_i} = \frac{1}{2}\frac{  A^{12}_{x_i}-A_{u_m}^{12} }{x_i-u_m}
+ \frac{1}{2} \sum_{\substack{k=1\\k\neq m}}^{g} \frac{ A^{12}_{u_k} - A_{u_m}^{12} }{u_k-u_m} \frac{\partial u_k}{\partial x_i}
-\frac{1}{2}\frac{\partial u_m}{\partial x_i} \!\!\sum_{\substack{a_j\in B \\a_j\neq u_m}}\frac{  \Aaj^{12}-\Aum^{12}}{a_j-u_m}\,, \quad1\leq m\leq g\,.\qquad
\label{Aum12}
\end{eqnarray}

Let us first prove \eqref{Aaj12} for $\Aaj$ with $a_j\notin\{u_1, \dots, u_{g}\}$.
Differentiating $\Aaj^{12}$ defined by \eqref{A12_aj_g} using Proposition \ref{proposition_Omega_aj} rewritten using $\Aaj^{12}$ as follows
\begin{equation*}
\frac{\partial \Omega(\paj)}{\partial x_i }  = \frac{\Omega(\paj)}{2(x_i-a_j)}\frac{\Axi^{12}}{\Aaj^{12}} \prod_{\alpha=1}^{g}\frac{x_i-u_\alpha}{a_j-u_\alpha}
\end{equation*}
 and \eqref{der_phi_aj_g}, we have
\begin{multline}
\label{Aaj_derivative}
\frac{\partial \Aaj^{12}}{\partial x_i} = \frac{\partial }{\partial x_i}\left\{ \frac{t}{4}\Omega(\paj) \varphi(\paj)\right\}
=\Aaj^{12}
\left\{  \frac{\partial_{x_i} \Omega(\paj)}{\Omega(\paj)} + \frac{\partial_{x_i} \varphi(\paj)}{\varphi(\paj)}
\right\}
\\
=\frac{\Aaj^{12}}{2}
\left\{  \frac{1}{x_i-a_j}\frac{\Axi^{12}}{\Aaj^{12}} \frac{\prod_{\alpha=1}^{g}(x_i-u_\alpha)}{\prod_{\alpha=1}^{g}(a_j-u_\alpha)}
 + \left( \frac{1}{a_j-x_i} + \sum_{\alpha=1}^{g}\frac{1}{a_j-u_\alpha} \frac{\partial u_\alpha}{\partial x_i} \right) \right\}.
\end{multline}
Using this in the left hand side of \eqref{Aaj12}, we see that it remains to prove the following equality:
\begin{equation*}
 \frac{\prod_{\alpha=1}^{g}(x_i-u_\alpha)}{\prod_{\alpha=1}^{g}(a_j-u_\alpha)} =1 + \frac{x_i-a_j}{\Axi^{12}} \sum_{k=1}^{g} \frac{ A^{12}_{u_k} }{u_k-a_j} \frac{\partial u_k}{\partial x_i}\,.
\end{equation*}
This is done by plugging in \eqref{umder_A} for the derivatives of $u_k$ and then using rational identity \eqref{rat1}.

Let us now prove \eqref{Aum12} with $1\leq m\leq g$. Differentiating $\Aaj^{12}$ defined by \eqref{A12_aj_g} with $a_j=u_m$ and using Proposition \ref{prop_Omega_um} rewritten using the notation of $\Aaj^{12}$,
 and \eqref{der_phi_um} for derivative of $\phi(\pum)$, we have
\begin{multline}
\label{Aum_derivative}
\frac{\partial \Aum^{12}}{\partial x_i} =
\frac{\Aum^{12}}{2}\left\{
-\frac{\partial u_m}{\partial x_i}
 \left( \frac{1}{x_i-u_m} - \sum_{\substack{ \alpha=1\\ \alpha\neq m}}^{g}\frac{1}{u_m-u_\alpha}  +\frac{1}{A_{u_m}^{12}}\sum_{\substack{a_j\in B\\a_j\neq u_m}} \frac{A_{a_j}^{12}}{a_j-u_m}     \right) \right.
\\
\left.
+\frac{1}{u_m-x_i} -\frac{\partial u_m}{\partial x_i}\sum_{\substack{a_j\in B\\a_j\neq u_m}} \frac{1}{u_m-a_j}  +\sum_{\substack{\alpha=1\\ \alpha \neq m}}^{g} \frac{1}{u_m-u_\alpha} \frac{\partial u_\alpha}{\partial x_i}
\right\}.
\end{multline}
Using this in the left hand side of \eqref{Aum12}, one sees that the terms in the second line and the sum over $a_j$ in the first line of \eqref{Aum_derivative} cancel against the corresponding terms in the right hand side of \eqref{Aum12}, and what remains to prove is the following identity
\begin{equation*}
\Aum^{12} \frac{\partial u_m}{\partial x_i}
 \left(  \sum_{\substack{ \alpha=1\\ \alpha\neq m}}^{g}\frac{1}{u_m-u_\alpha}      - \frac{1}{x_i-u_m} \right)
 =
\frac{  A^{12}_{x_i} }{x_i-u_m}
+ \sum_{\substack{k=1\\k\neq m}}^{g} \frac{ A^{12}_{u_k} }{u_k-u_m} \frac{\partial u_k}{\partial x_i}\,,
\end{equation*}
which is easily verified using \eqref{umder_A} for the derivatives of $u_k$ and rational identity \eqref{rat2} with $x=x_i$  and $N=g$ in the right hand side.

\bigskip

{\bf Acknowledgements.}  Although neither of the co-authors of this paper had Igor Moiseevich Krichever as their PhD advisor, both V.D. and V.S. have been strongly influenced by Krichever's works since their PhD times. This paper is yet another  example of that. We thank  Alexey Yu. Morozov and Patrick Labelle for very helpful discussions. V.D. acknowledges with gratitude the Simons Foundation grant no. 854861 and the Science Fund of Serbia grant IntegraRS. V.S. gratefully acknowledges
support from the Natural Sciences and Engineering Research Council of Canada through a Discovery grant and from the University of Sherbrooke.


\begin{thebibliography}{99}

\addcontentsline{toc}{section}{Bibliography}

	\bibitem{AhiezerAPPROX} N. I. Akhiezer, {\it Lekcii po Teorii Approksimacii}, [in Russian], OGIZ, Moscow-Leningrad, 1947, pp. 323.

\bibitem{Akh4} N. I. Akhiezer,  {\it Elements of the theory of elliptic functions},
			Translations of Mathematical Monographs,
			volume 79, American Mathematical Society, Providence, RI, 1990, pp. viii+237.


\bibitem{Apt1986}
A. I. Aptekarev, Asymptotic properties of polynomials orthogonal on a system of contours,
and periodic motions of Toda chains, Math.USSR. Sb 53 (1986), 233-260.


\bibitem{AF1995}
P. Argyres and A. Faraggi, {\it The vacuum structure and spectrum of N=2 supersymmetric
SU(N) gauge theory}, Phys. Rev. Lett. 73 (1995) 3931, hep-th/9411057


\bibitem{Bogatyrev2012} A. Bogatyrev, {\it Extremal polynomials and Riemann surfaces}, Springer Monographs in Mathematics, Springer, Heidelberg, 2012, pp.  xxvi+150 pp.


\bibitem{BDK}
S. A. Brazovskil, N. E. Dzyaloshinskil, and I. M. Krichever, {\it
Discrete Peierls models with exact solutions}, Zh. Eksp. Teor. Fiz. 83, 389-415; Sov. Phys. JETP 56(1), July 1982, 212-225.



\bibitem{DaTa}
E. Date, and S. Takana, {\it Analogue of inverse scattering theory for the discrete Hill's equation
and exact solutions for the periodic Toda lattice}, Progr. Theor. Phys. 55, 1976, pp. 447-465.


\bibitem{DHPh} E. D'Hoker, D. H.  Phong,
 {\it Lectures on supersymmetric Yang-Mills theory and integrable systems}, CRM Ser. Math. Phys.
Springer-Verlag, New York, 2002, 1-125.


\bibitem{DS2019} V. Dragovi\'{c}, V. Shramchenko, \textit{Algebro-geometric approach to an Okamoto transformation, the Painlev\'e VI and Schlesinger equations}, Annales Henri Poincar\'e, 2019,  20, no. 4, 1121-1148.

    \bibitem{isoharmonic} V. Dragovi\'c, V. Shramchenko, {\it Isoharmonic deformations and constrained Schlesinger systems}. {arXiv:2112.04110}

\bibitem{DR2019}  V. Dragovi\'c, M. Radnovi\'c, {\it Periodic Ellipsoidal Billiard Trajectories and Extremal
Polynomials}, Comm. Math. Phys. , 2019, Vol. 372, p. 183-211.

\bibitem{DR2023}  V. Dragovi\'c, M. Radnovi\'c, {\it Resonance of ellipsoidal billiard trajectories and extremal rational functions},  Advances in Mathematics, Article 109044, 51 p. Vol. 424, 2023.

    \bibitem{DMN} B. A. Dubrovin, V. B. Matveev, S. P. Novikov, \textit{Nonlinear equations of Korteweg-de Vries type, finite zone operators, and Abelian varieties}, Russ. Math. Surv. 31, No. 1, 1976. 59-146.

\bibitem{Dub1981} B. A. Dubrovin, \textit{Theta functions and non-linear equations}, Russ. Math. Surv. 1981,  36, No. 2 11-92.

    \bibitem{Fay92} J. Fay, {\it Kernel functions, analytic torsion, and moduli spaces}, Memoirs of the AMS \textbf{96} (1992), no 464.

  \bibitem{Fla1}  H. Flaschka,  {\it The Toda lattice: I. Existence of integrals}, Phys. Rev. B 9 1974 1924-1925
  \bibitem{Fla2} H.  Flaschka, {\it The Toda lattice: II}, Prog. Theor. Phys. 51, 1974, pp. 703-716.

\bibitem{FlaMcL} H. Flaschka H and D. McLaughlin, {\it  Canonically conjugate variables for the Korteweg -- de Vries equation and
the Toda lattice with periodic boundary conditions}, Prog. Theor. Phys. 55, 1976, 438-56.

 \bibitem{GKMMM}   A. Gorskii, I.M. Krichever, A. Marshakov, A. Mironov and A. Morozov, {\it Integrability
and Seiberg-Witten exact solution}, Phys. Lett. B355 (1995) 466, hep-th/9505035;

\bibitem{GS} P. G. Grinevich, M. U. Schmidt, {\it Period preserving nonisospectral flows and the moduli space of periodic solutions of soliton equations: The nonlinear Schr\"odinger equation}, Physica D, 87:73-98, 1995.
\bibitem{HK}
A. Henrici and T. Kappeler, {\it Global Birkhoff coordinates for the periodic Toda
lattice}, Nonlinearity {\bf 21} (2008), no.12, 2731-2758.

    \bibitem{Ince} E. L. Ince, {\it Ordinary differential equations} Dover, New York (1956).

    \bibitem{Japan} K. Iwasaki, H. Kimura, S. Shimomura, M. Yoshida, {\it From
Gauss to Painlev\'e. A modern theory of special functions}, Aspects
of Mathematics, Braunschweig (1991).

  \bibitem{KLTY}   A. Klemm, W. Lerche, S. Theisen, and S. Yankielowicz, {\it Simple singularities and
$N=2$ supersymmetric gauge theories}, Phys. Lett. B 344 (1995) 169, hep-th/9411058


\bibitem{KokoKoro} A. Kokotov, D. Korotkin, {\it A new hierarchy of integrable systems associated to Hurwitz spaces}, Philos. Trans. R. Soc. Lond. Ser. A Math. Phys. Eng. Sci. {\bf 366} (2008), no. 1867, 1055-1088.

    \bibitem{Krich1978} I. M. Krichever. \textit{Algebraic curves and non-linear difference equations}, Russ. Math. Surv. 1978, 33, No. 4 255-256.

\bibitem{Krich1981} I. M. Krichever,  \textit{The periodic non-Abelian Toda chain and its two-dimensional generalization}, pp. 82-89, Appendix in B. A. Dubrovin, \textit{Theta functions and nonlinear equations}, Russ. Math. Surv. 1981,  36, No. 2.


\bibitem{HillToda} H. P. McKean, P. van Moerbeke, {\it Hill and Toda curves},
Comm. Pure Appl. Math. 33 (1980), no. 1, 23-42.

\bibitem{Moer}
 P. van Moerbeke, {\it The spectrum of Jacobi matrices}, Inv. Math. 37. 1976, pp. 45-81.

\bibitem{NaTa}
T. Nakatsu and K.
Takasaki, {\it Whitham-Toda Hierarchy and $N = 2$ supersymmetric Yang-Mills Theory},
Mod. Phys. Lett. A11 (1996) 157, hep-th/9509162;

    \bibitem{Nov}
S. P. Novikov, S. V.  Manakov, L. P.  Pitaevskii, V. E. Zakharov,  {\it Theory of solitons},
Contemp. Soviet Math.,
Consultants Bureau [Plenum], New York, 1984, xi+276 pp.

\bibitem{PS1999}
F. Peherstorfer, K. Schiefermayr, {\it Description of extremal polynomials on several intervals
and their computation}, I, II, Acta Math. Hungar. 83 (1999), no. 1-2, 27-58, 59-83.


\bibitem{Per} L. V. Perkolab, {\it Inverse problem for periodic Jacobi matrices}, Theory of functions, functional analysis and their applications, (1984), Vol. 42, p. 107-121.

 \bibitem{SW1}   N. Seiberg and E. Witten, {\it Electro-magnetic duality, monopole condensation, and
confinement in $N = 2$ supersymmetric Yang-Mills theory}, Nucl. Phys. B426 (1994)
19-53, hep-th/9407087.

\bibitem{SW2}   N. Seiberg and E. Witten, {\it Monopoles, duality, and chiral symmetry breaking in $N = 2$
supersymmetric QCD}, Nucl. Phys. B431 (1994) 494, hep-th/9410167.


 \bibitem{Si2011} B. Simon, {\it Szeg\"o's Theorem and Its Descendants,} Princeton University Press, 2011.
\bibitem{Si2015a} B. Simon, {\it Harmonic Analysis. A Comprehensive Course in Analysis, Part 3}, AMS, 2015.
\bibitem{Si2015} B. Simon, {\it Operator Theory. A Comprehensive Course in Analysis, Part 4}, AMS, 2015.


 \bibitem{SY1992} M. L. Sodin, P. M. Yuditskii, {\it  Functions that deviate least from zero on closed subsets of the real axis}, Algebra i Analiz 4 (1992), no. 2, 1-61 (Russian, with Russian summary); English transl., St. Petersburg Math. J. 4 (1993), no. 2, 201-249.


\bibitem{To}
 M. Toda, {\it Studies of a non-linear lattice},
Phys. Rep. 18C (1975), no. 1, 1-123.


\bibitem{Widom} H. Widom, {\it Extremal polynomials associated with a system of curves in the complex plane}, Adv. Math. 3, (1969), 127-232.


\end{thebibliography}
\end{document}